\documentclass[12pt]{amsart}
\usepackage{amsfonts}
\usepackage{amssymb}
\numberwithin{equation}{section}
\theoremstyle{plain}
 \newtheorem{thm}{Theorem}[section]
 \newtheorem{lem}[thm]{Lemma}
 \newtheorem{cor}[thm]{Corollary}
 \newtheorem{prop}[thm]{Proposition}
 
\theoremstyle{definition}
\newtheorem{defn}[thm]{Definition}
 
 \newtheorem{rem}[thm]{Remark}

\newcommand{\al}{\alpha}
\newcommand{\bt}{\beta}
\newcommand{\gm}{\gamma}

\newcommand{\ld}{\lambda}

\newcommand{\Ph}{\Phi}

\newcommand{\Up}{\Upsilon}

\newcommand{\ps}{\psi}

\newcommand{\q}{\quad}

\newcommand{\wh}{\widehat}
\newcommand{\wt}{\widetilde}

\newcommand{\R}{\mathbb{R}}

\newcommand{\rd}{\mathbb R ^d}

\newcommand{\law}{\mathcal L}
\newcommand{\sek}{\int_0^{\infty}}
\newcommand{\pn}{\par\noindent}

\newcommand{\mM}{\mathcal{M}}
\newcommand{\mE}{\mathcal{E}}
\newcommand{\mN}{\mathcal{N}_{\alpha}}
\newcommand{\n}{\noindent}
\newcommand{\bR}{\mathbb{R}}

\begin{document}
\setlength{\baselineskip}{18pt}
\setlength{\parindent}{1.8pc}

\n
{\bf \large New classes of infinitely divisible distributions related to the Goldie--Steutel--Bondesson class }

\vskip 5mm
\n
Takahiro Aoyama $\cdot$ Alexander Lindner $\cdot$ Makoto Maejima\\

\addtocounter{footnote}{1}\footnotetext{Takahiro Aoyama\\
Department of Mathematics,
Tokyo University of Science, 2641, Yamazaki, Noda 278-8510, Japan\\
e-mail: aoyama$\_{}$takahiro@ma.noda.tus.ac.jp}
\addtocounter{footnote}{1}\footnotetext{Alexander Lindner\\
Technische Universit\"at
Braunschweig, Institut f\"ur Mathematische Stochastik,
Pockelsstra{\ss}e 14, D-38106 Braunschweig, Germany\\
e-mail: a.lindner@tu-bs.de}

\addtocounter{footnote}{1}\footnotetext{Makoto Maejima (Corresponding author)\\
Department of Mathematics,
Keio University, 3-14-1, Hiyoshi, Kohoku-ku, Yokohama 223-8522, Japan\\
e-mail: maejima@math.keio.ac.jp}
\overfullrule=0pt

\vskip 10mm
\n
(Running head: New classes of infinitely divisible distributions)

\vskip 10mm
\n
{\bf Abstract}
Recently, many classes of infinitely divisible distributions on
$\rd$ have been characterized in several ways. Among others, the
first way is to use L\'evy measures, the second one is to use
transformations of L\'evy measures, and the third one is to use
mappings of infinitely divisible distributions defined by stochastic
integrals with respect to L\'evy processes. In this paper, we are
concerned with a class of mappings, by which we construct new
classes of infinitely divisible distributions on $\rd$. Then we
study a special case in $\R^1$, which is the class of infinitely
divisible distributions without Gaussian parts generated by
stochastic integrals with respect to a fixed compound Poisson
processes on $\R^1$. This is closely related to the
Goldie--Steutel--Bondesson class.

\vskip 5mm
\n

\allowdisplaybreaks
\vskip 10mm
%%%%%%%%%%%%%%%%% Section 1 %%%%%%%%%%%%%%%%%
\section{Introduction}

Throughout this paper, $\law (X)$ denotes the law of an $\rd$-valued random variable $X$ and
$\wh \mu (z), z\in\rd$, denotes the characteristic function of a probability distribution $\mu$ on $\rd$.
Also $I(\rd)$ denotes the class of all infinitely divisible distributions on $\rd$,
$I_{{\rm sym}}(\rd) =\{ \mu\in I(\rd) :\mu\,\,\text{is symmetric on}\,\, \rd\} $,
$I_{\log}(\rd)= \{ \mu\in I(\rd) :\int_{|x|>1}\log |x|\mu(dx)<\infty\} $ and
$I_{\log^m}(\rd)= \{ \mu\in I(\rd) :\int_{|x|>1}(\log |x|)^m\mu(dx)<\infty\} $,
where $|x|$ is the Euclidean norm of $x\in \rd$.
Let $C_{\mu}(z), z\in \rd$, be the cumulant function of $\mu\in I(\rd)$.
That is, $C_{\mu}(z)$ is a continuous function
with $C_{\mu}(0)=0$ such that $\wh \mu(z) = \exp \left \{ C_{\mu}(z)\right \}, z\in \rd$.

We use the generating triplet $(A, \nu, \gm)$ of $\mu\in I(\rd)$ in the sense that
\begin{align*}
C_{\mu}(z) = -{2^{-1}} \langle z,Az \rangle & + {\rm i}
\langle \gm,z \rangle\\
&+ \int_{\rd}\left(e^{{\rm i} \langle z,x \rangle }-1-
{{\rm i} \langle z,x \rangle}(1+|x|^2)^{-1}\right)\nu(dx),\,z\in\rd, \nonumber
\end{align*}
where $A$ is a symmetric nonnegative-definite $d \times d$ matrix,
$\gamma\in\rd$ and $\nu$ is a measure  (called the L\'evy measure) on
$\rd$ satisfying
\begin{equation*}
\nu(\{0\}) = 0 \,\,\text{and}\,\,\int_{\rd} (|x|^{2} \wedge 1) \nu(dx) < \infty.
\end{equation*}

The polar decomposition of the L\'evy measure $\nu$ of $\mu\in I(\rd)$, with $0<\nu(\R^d)\le\infty$,
is the following:
There exist a measure $\ld$ on $S=\{\xi\in\rd : |\xi |=1\}$
with $0<\ld(S)\le\infty$ and
a family $\{\nu_{\xi}\colon \xi\in S\}$ of measures on $(0,\infty)$ such that
$\nu_{\xi}(B)$ is measurable in $\xi$ for each $B\in\mathcal B((0,\infty))$,
$0<\nu_{\xi}((0,\infty))\le\infty$ for each $\xi\in S$,
\begin{align} \label{eq-2.1}
\nu(B)=\int_S \ld(d\xi)\int_0^{\infty} 1_B(r\xi)\nu_{\xi}(dr), \,\,
B\in \mathcal B (\mathbb R^d \setminus \{ 0\}).
\end{align}
Here $\ld$ and $\{\nu_{\xi}\}$ are uniquely determined by $\nu$ up
to multiplication by a measurable function $c(\xi)$ and
${c(\xi)}^{-1}$ with $0<c(\xi)<\infty$. The measure $\nu_\xi$ is  a
L\'evy measure on $(0,\infty)$ for $\lambda$-a.e. $\xi\in S$. We say
that $\nu$ has the polar decomposition $(\ld ,\nu_{\xi})$ and
$\nu_{\xi}$ is called the radial component of $\nu$. (See, e.g.,
Lemma 2.1 of \cite{BNMS06} and its proof.)

The classes which we are going to study in this paper are the following.

\begin{defn}
(Class $E_{\al}(\rd), \al >0.)$ We say that $\mu\in
I(\rd)$ belongs to the class $E_{\al}(\rd)$ if $\nu =0$ or $\nu
\ne 0$ and, in case $\nu\ne 0$, $\nu_{\xi}$ in \eqref{eq-2.1}
satisfies
\begin{equation*}
\nu_{\xi}(dr) = r^{\al -1}g_{\xi}(r^{\al})dr , \,\, r>0,
\end{equation*}
for some function $g_{\xi}(r)$, which is completely monotone in $r\in (0,\infty)$ for
$\ld$-a.e.~$\xi$, is measurable in $\xi$ for each $r>0$ and satisfies
$$
\int _0^{\infty} (r^{\al+1} \wedge r^{\al-1})g_{\xi}(r^{\al})dr < \infty,\q r>0,\q \text{$\ld$-a.e. $\xi$}.
$$
\end{defn}

The following four known classes are needed in our discussion.

\vskip 3mm \n (1) Class $B(\rd) $ (the Goldie--Steutel--Bondesson
class): $\mu\in I(\rd)$ belongs to the class $B(\rd)$ if $\nu=0$ or
$\nu\ne 0$ and, in this case, $\nu_{\xi}$ in \eqref{eq-2.1}
satisfies
\begin{equation*}
\nu_{\xi}(dr)=g_{\xi}(r)dr,
\end{equation*}
where $g_{\xi}(r)$ is completely monotone in $r\in (0,\infty)$ for
$\ld$-a.e.~$\xi$ and is measurable in $\xi$ for each $r>0$. Hence
$E_1(\rd)= B(\rd)$.

\vskip 3mm \n (2) Class $L(\rd)$ (the class of selfdecomposable
distributions): $\mu\in I(\rd)$ belongs to the class $L(\rd)$ if
$\nu=0$ or $\nu\ne 0$ and, in this case, $\nu_{\xi}$ in
\eqref{eq-2.1} satisfies
\begin{equation*}
\nu_{\xi}(dr)=r^{-1}{k_{\xi}(r)}dr,
\end{equation*}
where $k_{\xi}(r)$ is nonincreasing in $r\in (0,\infty)$ for
$\ld$-a.e.~$\xi$ and is measurable in $\xi$ for each $r>0$. \vskip
3mm \n (3) Class $M(\rd)$: $\mu\in I(\rd)$ belongs to the class
$M(\rd)$ if $\nu=0$ or $\nu\ne 0$ and, in this case, $\nu_{\xi}$ in
\eqref{eq-2.1} satisfies
\begin{equation*}
\nu_{\xi}(dr)= r^{-1}{g_{\xi}(r^2)}dr,
\end{equation*}
where $g_{\xi}(r)$ is completely monotone in $r\in (0,\infty)$ for
$\ld$-a.e.~$\xi$ and is measurable in $\xi$ for each $r>0$.

\vskip 3mm \n (4) Class $T(\rd)$ (the Thorin class): $\mu\in I(\rd)$
belongs to the class $T(\rd)$ if $\nu=0$ or $\nu\ne 0$ and, in this
case, $\nu_{\xi}$ in \eqref{eq-2.1} satisfies
\begin{equation*}
\nu_{\xi}(dr)=r^{-1}{g_{\xi}(r)}dr,
\end{equation*}
where $g_{\xi}(r)$ is completely monotone in $r\in (0,\infty)$ for
$\ld$-a.e.~$\xi$ and is measurable in $\xi$ for each $r>0$.

\vskip 3mm

We introduce four mappings from $I(\rd)$ (or $I_{\log}(\rd)$) into
$I(\rd)$, which are related to the classes above. Throughout this
paper,  $\{X_t^{(\mu)}\}$ denotes a L\'evy process on $\rd$ with
$\law (X_1^{(\mu)}) =\mu$.

\begin{defn}
(1) For $\al >0$,
$
\mE _{\al}(\mu)=\law\left (\int_0^1(\log t^{-1})^{1/\al}dX_t^{(\mu)}\right ),\,\,\mu\in I(\rd).
$

\n
(2)
$
\Phi (\mu) = \law\left (\sek e^{-t}dX_t^{(\mu)}\right ),\,\,\mu\in I_{\log}(\rd).
$

\n
(3)
$
\mM (\mu) = \law\left (\int_0^{\infty}m^*(t) dX_t^{(\mu)}\right ),\,\,\mu\in I_{\log}(\rd),
$
where $m(x) = \int_x^{\infty}u^{-1}e^{-u^2}du, x>0$, and $m^*(t)$ is its inverse function in the sense that
$m(x)=t$ if and only if $x=m^*(t)$.

\n
(4)
$
\Psi (\mu) =\law\left ( \int _0^{\infty} e^*(t) dX_t^{(\mu)}\right ), \,\, \mu\in I_{\log}(\rd),
$
where $e(x) = \int_x^{\infty}u^{-1}{e^{-u}}du, x>0,$ and  $e^*(t)$ is its inverse function in the sense that
$e(x)=t$ if and only if $x=e^*(t)$.
\end{defn}
Only the mapping $\mE_{\al}$ (for $\alpha\neq 1$) is new.
%(See Theorem 3.1 below.)
It is known that  $\frak D (\Phi) = \frak D (\mathcal M ) = \frak D
(\Psi) = I_{\log}(\rd)$, where the domain $\frak D (*)$ means the
set of infinitely divisible distributions $\mu$ on $\mathbb{R}^d$ on
which the $*$-mapping is definable, in the sense of improper
integrals with respect to independently scattered random measures on
$\mathbb{R}^d$, as
in Definitions~2.3 and~3.1 of Sato~\cite{S07}. (For the
determination of $\frak D(\Phi)$, $\frak D(\Psi)$ and $\frak
D(\mM)$, see \cite{W82}, \cite{BNMS06} and \cite{AMR08},
respectively.) For $\mE _{\al}$, we have $\frak D (\mE _{\al}) =
I(\rd)$, as shown in Proposition~\ref{prop-1} below.

\begin{rem} \label{rem-1.3}
$\mE_{1}$ is known as the Upsilon mapping (denoted by $\Up$ in the literature)
and it is known that $\frak D(\Up) =I(\rd)$ and
$\Up(I(\rd)) = B(\rd)$.
Recall that $E_1(\rd) = B(\rd)$.
Hence
\begin{equation} \label{eq-1.3}
E_1(\rd) = \mE_1(I(\rd)).
\end{equation}
\end{rem}

The paper is organized as follows. In Section 2, we show several
properties of the mapping $\mE_{\al}$. In Section 3, we show that
$E_{\al}(\rd) = \mE_{\al}(I(\rd)),\,\, \al >0$. This gives
us stochastic integral representations of the elements of the class
$E_{\al}(\rd)$. In Section~4, we consider the composition $\mN$, say, of two
mappings $\Phi$ and $\mE_{\al}$, and show in particular that $\mM(=\mathcal{N}_2)$ is the composition of $\Phi$ and
$\mE_2$. Then as an application of this equality, we show that
the limit of certain subclasses of $N_{\al}(\rd) := \mN(I_{\log}(\rd))$, constructed by the iteration of the
mapping $\mN$, is
the closure of the class of the stable distributions as Maejima and Sato \cite{MS08}
showed for other mappings. In Section 5, we restrict
ourselves to the case $d=1$ and characterize
\begin{align}
E_{\al}^0 (\bR^1) & :=   \{ \mu \in E_\alpha (\mathbb{R}^1) :
\mbox{$\mu$ has no Gaussian part}\}
\label{def-0}\\
E_{\al}^{0, \rm sym}(\mathbb R^1) & :=   E_\alpha^0 (\mathbb{R}^1)
\cap I_{\rm sym}(\mathbb R ^1) \label{def-sym}
\end{align}
and certain subclasses of $E_\alpha^0(\bR^1)$ which correspond to
L\'evy processes of bounded variation with zero drift,
by (essential improper) stochastic integrals
with respect to some compound Poisson processes.
This gives us a new
sight of the Goldie--Steutel--Bondesson class in $\mathbb{R}^1$.

\vskip 10mm
%%%%%%%%%%%%%%%%% Section 2 %%%%%%%%%%%%%%%%%
\section{Several properties of the mapping $\mE_{\al}$}

We start with showing several properties of the mapping $\mE_{\al}$.

\begin{prop} \label{prop-1}
Let $\alpha >0$. \\
$(i)$ $\mE_{\al}(\mu)$ can be defined for any $\mu \in I(\R^d)$ and is
infinitely divisible, and we have \linebreak $\int_0^1 |C_\mu (z
(\log t^{-1})^{1/\alpha})|\, dt < \infty$ and
\begin{equation*}
C_{\mE_\alpha (\mu) }(z) = \int_0^1 C_\mu ( z \left(\log
t^{-1})^{1/\alpha}\right) \, dt, \quad z \in \R^d.
\end{equation*}
$(ii)$ The generating triplet $(\widetilde{A}, \widetilde{\nu},
\widetilde{\gamma})$ of $\wt\mu=\mE_\alpha (\mu)$ can be calculated from
$(A,\nu,\gamma)$  of $\mu$ by
\begin{align}
\widetilde{A} & = \Gamma (1 + 2/\alpha) \, A ,
%\label{eq-tripel-A}
\nonumber\\
 \widetilde{\nu} (B) & =  \int_0^\infty \nu(
u^{-1} B) \alpha u^{\alpha-1} e^{-u^\alpha} \, du, \quad B \in
\mathcal{B} (\R^d\setminus\{ 0\}),
\label{eq-tripel-nu}\\
\wt{\gamma} & = \Gamma (1 + 1/\alpha) \, \gamma + \int_0^\infty
\alpha u^\alpha e^{-u^\alpha} \int_{\R^d} x \left(
\frac{1}{1+|ux|^2} - \frac{1}{1+|x|^2} \right) \, \nu(dx) \, du .
\label{eq-tripel-gamma}
\end{align}
$(iii)$ The mapping $\mE_\alpha : I(\R^d) \to I(\R^d)$ is
one-to-one.\\
$(iv)$ Let $\mu_n \in I(\R^d)$, $n=1,2, \ldots$ If $\mu_n$ converges
weakly to $\mu \in I(\R^d)$ as $n\to\infty$, then $\mE_\alpha
(\mu_n)$ converges weakly to $\mE_\alpha (\mu)$ as $n\to\infty$.
Conversely, if $\mE_\alpha (\mu_n)$ converges weakly to $\wt{\mu}$
for some distribution $\wt{\mu}$ as $n\to\infty$, then $\wt{\mu} =
\mE_\alpha (\mu)$ for some $\mu \in I(\R^d)$ and $\mu_n$ converges
weakly to $\mu$ as $n\to\infty$. In particular, the range
$\mE_\alpha(I(\R^d))$ is closed under weak convergence.\\
$(v)$ For any $\mu \in I(\R^d)$ we also have $$\mE_\alpha(\mu) =
\law \left( \int_0^1 \left(\log \frac{1}{1-t}\right)^{1/\alpha} \,
dX_t^{(\mu)} \right) = \law \left( \lim_{s\downarrow 0} \int_s^1
\frac{1}{\alpha t} (\log t^{-1})^{1/\alpha -1} X_t^{(\mu)} \, dt
\right),$$ where the limit is almost sure.
\end{prop}

\begin{proof}
(The proof follows along the lines of Proposition 2.4 of \cite{BNMS06}.
However, we give the proof for the completeness of the paper.)

(i) The function $f(t) = (\log t^{-1})^{1/\alpha}
\mathbf{1}_{(0,1]}(t)$ is clearly square integrable, hence the result
follows from Sato~\cite{S06}, see also Lemma~2.3 in
Maejima~\cite{M07}.

(ii) By a general result (see Lemma~2.7 and Corollary~4.4 of
Sato~\cite{S04}) and a change of variable,  we have
$$\wt{A} = \left ( \int_0^1 (\log t^{-1})^{2/\alpha} \, dt\right  ) \, A =
\left ( \int_0^\infty u^{2/\alpha} e^{-u} \, du\right ) \, A = \Gamma ( 1 +
2/\alpha) \, A,$$
\begin{align*}
\wt{\nu}(B) = \int_0^1 \nu((\log t^{-1})^{-1/\al}B)dt
= \int_0^\infty \nu(u^{-1} B) \, \al u^{\al - 1}e^{-u^\al}du,
\end{align*}
\begin{align*}
\wt{\gm}&=\int_0^1 (\log t^{-1})^{1/\al} \left(\gm + \int_{\rd} x
\left( \frac{1}{1+|(\log t^{-1})^{1/\alpha} x|^2} - \frac{1}{1+|x|^2} \right)\nu(dx) \right)dt \\
&=\gm \int_0^\infty v^{1/\al} e^{-v} \, dv + \int_0^\infty \al u^\al e^{-u^\al} \int_{\rd}
\left(\frac{1}{1+|ux|^2} - \frac{1}{1+|x|^2} \right)\nu(dx)\,du.
\end{align*}

(iii) By (i), we have for each $z\in\R^d$,
$$C_{\mE_\alpha (\mu)} (z) = \int_0^1 C_\mu (z (\log
t^{-1})^{1/\alpha})\, dt = \int_0^\infty C_\mu (z
v^{1/\alpha})e^{-v} \, dv.$$ Hence we conclude that for each $u >
0$ and $z \in \R^d$,
$$\frac{1}{u} C_{\mE_\alpha (\mu)} (u^{-1/\alpha} z) =
\int_0^\infty \frac{1}{u} C_\mu \left( \left( \frac{v}{u}
\right)^{1/\alpha} z\right) e^{-v} \, dv = \int_0^\infty C_\mu
(w^{1/\alpha} z) e^{-uw} \, dw.$$ Hence we see that for each
$z\in\R^d$, the function $(0,\infty) \to \R$, $u \mapsto u^{-1}
C_{\mE_\alpha (\mu)}(u^{-1/\alpha} z)$ is the Laplace transform of
$(0,\infty) \to \R$, $w \mapsto C_\mu (w^{1/\alpha} z)$. Hence for
each fixed $z\in\R^d$, $C_\mu(w^{1/\alpha} z)$ is determined by
$\mE_\alpha (\mu)$ for almost every $w\in (0,\infty)$, and by
continuity for every $w > 0$. In particular for $w=1$, we see that
$C_\mu(z)$ is determined by $\mE_\alpha (\mu)$ for every $z\in\R^d$.

(iv) Apart from minor adjustments, the proof is the same as that
of Proposition~2.4 (v) in Barndorff-Nielsen et al.\,\cite{BNMS06} and
hence omitted.

(v) The first equality is clear by duality (see Sato~\cite{S99},
Proposition~41.8). For the second, observe that $\int_s^1 (\log
t^{-1})^{1/\alpha} \, dX_t^{(\mu)}$ converges almost surely to
$\int_0^1 (\log t^{-1})^{1/\alpha} \, dX_t^{(\mu)}$ as $s\downarrow
0$ by the independently scattered random measure property of
$X_t^{(\mu)}$. Using partial integration, we conclude
$$\int_s^1 (\log t^{-1})^{1/\alpha} \, dX_t^{(\mu)} = -
X_s^{(\mu)} (\log s^{-1})^{1/\alpha} - \int_s^ 1 X_t^{(\mu)} d(\log
t^{-1})^{1/\alpha}.$$ But $\lim_{s\downarrow 0} X_s^{(\mu)} (\log
s^{-1})^{1/\alpha} = 0$ a.s. (see Sato~\cite{S99},
Proposition~47.11), and the claim follows.
\end{proof}

\begin{cor} \label{Cor-2}
Let $\alpha> 0$. Then a distribution $\mu$ is symmetric if and only
if $\mE_\alpha (\mu)$ is symmetric.
\end{cor}

\begin{proof}
It is well known that an infinitely divisible distribution $\mu$
with generating triplet $(A,\nu,\gamma)$ is symmetric if and only if
$\nu$ is symmetric and $\gamma=0$. Symmetry of $\mu$ hence implies
symmetry of $\mE_\alpha(\mu)$ by \eqref{eq-tripel-nu} and
\eqref{eq-tripel-gamma}. Conversely, suppose that
$\widetilde{\mu}=\mE_\alpha(\mu)$ with triplet $(\widetilde{A},
\widetilde{\nu},\widetilde{\gamma})$ is symmetric. Then
$$\int_0^\infty \nu(u^{-1} B) \alpha u^{\alpha-1} e^{-u^\alpha}
\, du = \int_0^\infty \nu(-u^{-1} B) \alpha u^{\alpha-1}
e^{-u^\alpha} \, du$$ for every $B\in \mathcal{B}(\R^d\setminus \{
0 \})$. In particular,
$$\int_0^\infty \nu(u^{-1}t^{-1} B) u^{\alpha-1} e^{-u^\alpha} \,
du = \int_0^\infty \nu(-u^{-1}t^{-1} B) u^{\alpha-1} e^{-u^\alpha}
\, du \quad \forall \; t > 0,$$ so that
$$\int_0^\infty \nu(w^{-1} B) w^{\alpha-1} e^{-w^\alpha/t^\alpha} \,
dw = \int_0^\infty  \nu(-w^{-1} B) w^{\alpha-1}
e^{-w^\alpha/t^\alpha} \, dw \quad \forall\; t > 0.$$ By the
uniqueness theorem for the Laplace transform, it follows that for
fixed $B$, $\nu(w^{-1} B) = \nu(-w^{-1}B)$ for almost every $w>0$.
Now let $B$ be of the form
$$B = \{ x \in \R^d : |x|> r \quad \mbox{and}\quad x/|x| \in  U\}$$
for some $r>0$ and some $U\in \mathcal B (S)$. Then both $u
\mapsto \nu(u B)$ and $v\mapsto \nu(-u B)$ are c\`adl\`ag, and we
conclude equality of $\nu(B)$ and $\nu(-B)$. This shows that $\nu$
is symmetric, which then also shows $\gamma = \widetilde{\gamma} =
0$ by \eqref{eq-tripel-gamma}.
\end{proof}

\allowdisplaybreaks
\vskip 10mm
%%%%%%%%%%%%%%%%%%%%%%%%%%%% Section 3 %%%%%%%%%%%%%%%%%%%%%%%%%%%%
\section{Stochastic integral characterization of the classes $E _{\al} (\rd)$}

We start with stating the following known result. In what follows,
for two mappings $\Phi_1$ and $\Phi_2$, $\Phi_1\circ \Phi_2$ means
their composition $(\Phi_1\circ\Phi_2)(\mu) = \Phi_1(\Phi_2(\mu))$.

\begin{thm}
$(1)$
$L(\rd) = \Phi (I_{\log}(\rd))$. {\rm (\cite{W82} and others.)}

\n
$(2)$
$M(\rd) = \mM (I_{\log}(\rd))$.

\n
$(3)$
$\Phi \circ \mE _1 = \Psi$ and $T(\rd)=\Psi(I_{\log}(\rd))$. {\rm (\cite{BNMS06}.)}
\end{thm}

In \cite{AMR08}, $M(\rd)\cap I_{\rm sym} (\rd)$ is studied.
The statement (2) above can be shown by exactly the same way as in \cite{AMR08}.

Now we want to prove  the following two theorems.

\begin{thm}\label{theorem:th2}
For any $0<\al<\bt$,
\begin{align*}
E_{\al}(\rd)\subset E_{\bt}(\rd).
\end{align*}
\end{thm}

The following is an extension of \eqref{eq-1.3} in
Remark~\ref{rem-1.3} for general $\al >0$.
\begin{thm}\label{theorem:th1}
For $ \al >0$,
\begin{align*}
E_{\al}(\rd)=\mE _{\al}(I(\rd)).
\end{align*}
\end{thm}

\vskip 3mm \pn {\it Proof of Theorem \ref{theorem:th2}.} Let $0<
\alpha < \beta $. Then if $\mu\in E_{\al}(\rd)$, $\nu_{\xi}$
of $\mu$ is
\begin{align*}
\nu_{\xi}(dr)=r^{\al-1}g_{\xi}(r^{\al})dr
=r^{\bt-1}\frac{g_{\xi}\left((r^{\al/\bt})^{\bt}\right)}{r^{\bt-\al}}
=r^{\bt-1}\frac{g_{\xi}\left((r^{\al/\bt})^{\bt}\right)}
{\left(r^{(\bt-\al)/\bt}\right)^{\bt}}.
\end{align*}
Let
\begin{align*}
h_{\xi}(x)=\frac{g_{\xi}(x^{\al/\bt})}{x^{(\bt-\al)/\bt}}.
\end{align*}
Note that if $g$ is completely monotone and $\psi$ a nonnegative
function such that $\ps '$ is completely monotone, then the
composition $g\circ\ps $ is completely monotone (see, e.g.,
Feller~\cite{F66}, page 441, Corollary~2), and if $g$ and $f$ are
completely monotone then $gf$ is completely monotone. Thus
$g_{\xi}(x^{\al/\bt})$ is completely monotone and then $h_{\xi}(x)$
is also completely monotone, and we have
\begin{align*}
\nu_{\xi}(dr)=r^{\bt-1}h_{\xi}(r^{\bt}).
\end{align*}
Hence  $\mu\in E_{\bt}(\rd)$.
\qed

\vskip 3mm
\pn
{\it Proof of Theorem \ref{theorem:th1}.}

\n (i) (Proof for that $E_{\al}(\rd) \supset\mE_{\al}(I(\rd))$.) Let
$\wt{\mu}\in \mE _{\al}(I(\rd))$. Then $\wt\mu = $\linebreak $\law\left (
\int_0^1 (\log t^{-1})^{1/\al}dX_t^{(\mu)}\right )$ for some $\mu\in
I(\rd)$, and hence
\begin{align*}
\wt{\nu}(B):=\nu_{\wt{\mu}}(B)%=\int_0^1\nu
%\left((\log t^{-1})^{-1/\al}B\right)dt\\
=
\al\int_0^{\infty}\nu(u^{-1}B)u^{\al-1}e^{-u^{\al}}du,\label{nual}
\end{align*}
where $\nu$ is the L\'evy measure of $\mu$ and below $\nu_{\xi}$ is
the radial component  of $\nu$. Thus, the spherical component
$\widetilde{\lambda}$ of $\widetilde{\nu}$ is equal to the spherical
component $\lambda$ of $\nu$, and the radial component
$\wt\nu_{\xi}$ of $\wt\nu$ satisfies that, for
$B\in\mathcal{B}\left((0,\infty)\right)$
\begin{align}
\wt{\nu}_{\xi}(B)&=\al\int_0^{\infty}u^{\al-1}e^{-u^{\al}}du\int_0^{\infty}
\mathbf{1}_B(xu)\nu_{\xi}(dx)
%\label{eq-nu-wtnu}
\nonumber\\
%&=\al\int_0^{\infty}\nu_{\xi}(dx)\int_0^{\infty}1_B(xu)u^{\al-1}e^{-u^\al}du\nonumber\\
&=\al\int_0^{\infty}\nu_{\xi}(dx)\int_0^{\infty}\mathbf{1}_B(y)(y/x)^{\al-1}
e^{-(y/x)^{\al}}x^{-1}dy\nonumber\\
%&=\int_0^{\infty}1_B(y)y^{\al-1}dy\int_0^{\infty}\al x^{-\al}e^{-(y/x)^{\al}}\nu_{\xi}(dx)\nonumber\\
&=: \int_0^{\infty}\mathbf{1}_B(y)y^{\al-1}\wt
g_{\xi}(y^{\al})dy,\nonumber
\end{align}
where
\begin{align*}
\wt g_{\xi}(r)=\int_0^{\infty}\al
x^{-\al}e^{-r/x^{\al}}\nu_{\xi}(dx) =\int_0^{\infty}e^{-ru}\wt
Q_{\xi}(du),
\end{align*}
with the measure $\widetilde{Q}_\xi$ being defined by
\begin{align*}
\wt
Q_{\xi}(B)=\al\int_0^{\infty}\mathbf{1}_B(x^{-\al})x^{-\al}\nu_{\xi}(dx),
\quad B \in \mathcal{B}((0,\infty)).
\end{align*}
We conclude that $\wt g_{\xi}(\cdot)$ is completely monotone. Thus,
\begin{align*}
\wt{\nu}_{\xi}(dy)=y^{\al-1}\wt g_{\xi}(y^{\al})dy
\end{align*}
for some completely monotone function $\wt g_{\xi}$. This concludes
that $\wt\mu\in E_{\al}(\rd)$.

\n (ii) (Proof for that $E_{\al}(\rd) \subset\mE_{\al}(I(\rd))$.) Let
$\widetilde{\mu} \in E_{\al}({\rd})$ with L\'evy measure
$\widetilde{\nu}$ of the form
$$\wt{\nu} (B) = \int_{S} \wt{\lambda} (d\xi) \int_0^\infty \mathbf{1}_B
(r\xi) r^{\alpha-1} \wt g_\xi(r^\alpha) dr, \quad B\in \mathcal{B}
(\R^d\setminus \{ 0 \}),$$ where $g_\xi(r)$ is completely monotone
in $r$ and measurable in $\xi$. For each $\xi$, there exists a Borel
measure $\widetilde{Q}_\xi$ on $[0,\infty)$ such that $\wt g_\xi(r)
= \int_{[0,\infty)} e^{-rt} \, \widetilde{Q}_\xi(dt)$ and $\wt
Q_\xi(B)$ is measurable in $\xi$ for each $B \in
\mathcal{B}([0,\infty))$ (see the proof of Lemma~3.3 in
Sato~\cite{S80}). For $\wt{\nu}$ to be a L\'evy measure, it is
necessary and sufficient that
\begin{align*}
\infty & >  \int_S \wt{\lambda} (d\xi) \int_0^1 r^{\alpha + 1}
\wt g_\xi(r^\alpha) \, dr + \int_S \wt{\lambda} (d\xi) \int_1^\infty
r^{\alpha -1} \wt g_\xi(r^\alpha) \, dr \\
& =  \int_S \wt{\lambda} (d\xi) \int_0^1 r^{\alpha + 1}dr
\int_{[0,\infty)} e^{-r^\alpha t} \wt{Q}_\xi(dt) \,  \\
&  \hskip 20mm  + \int_S \wt{\lambda} (d\xi) \int_1^\infty r^{\alpha-1}dr
\int_{[0,\infty)} e^{-r^\alpha t} \wt{Q}_\xi(dt) \, \\
& =  \int_S \wt{\lambda} (d\xi) \alpha^{-1} \int_{[0,\infty)}
t^{-1-2/\alpha}\wt{Q}_\xi (dt)
\int_0^t u^{2/\alpha} e^{-u} \, du \, \\
& \hskip 20mm + \int_S \wt{\lambda} (d\xi) \alpha^{-1} \int_{[0,\infty)}
t^{-1} e^{-t} \wt{Q}_\xi(dt),
\end{align*}
where we have used Fubini's theorem and the substitution $u=r^\alpha t$
in the first term. From this it is easy to see that $\wt{\nu}$ is
a L\'evy measure if and only if $\int_S \wt{\lambda} (d\xi)
\wt{Q}_\xi ( \{ 0 \}) = 0$ (which we shall assume without comment
from now on) and
\begin{equation} \label{eq-Q2}
\int_S \wt{\lambda} (d\xi) \int_0^1 t^{-1} \wt{Q}_\xi (dt) <
\infty, \quad \int_S \wt{\lambda} (d\xi) \int_1^\infty
t^{-1-2/\alpha} \wt{Q}_\xi (dt) < \infty.
\end{equation}
In part (i) we have defined $\widetilde{Q}_\xi = U(\rho_\xi)$ as the
image measure of $\rho_\xi$ under the mapping $U:(0,\infty)\to
(0,\infty), r \mapsto r^{-\alpha}$, where $\rho_\xi$ has density
$r\mapsto \alpha r^{-\alpha}$ with respect to $\nu_\xi$. Denoting by
$V:r \mapsto r^{-1/\alpha}$, the inverse of $U$, it follows that
$\rho_\xi$ is the image measure of $\wt{Q}_\xi$ under the mapping
$V$. Hence, given $\wt{Q}_\xi$, we define $\nu_\xi$ as having
density $r\mapsto \alpha^{-1} r^{\alpha}$ with respect to the image
measure $V(\wt Q_\xi)$ of $\wt Q_\xi$ under $V$, i.e.
$$\nu_\xi(B) = \alpha^{-1} \int_0^\infty \mathbf{1}_B (r^{-1/\alpha})
r^{-1} \wt{Q}_\xi (dr), \quad B\in\mathcal{B} ((0,\infty)).$$
Define further a measure $\nu$ to have spherical component
$\lambda = \widetilde{\lambda}$ and radial parts $\nu_\xi$, i.e.
$$\nu(B) = \int_S \widetilde{\lambda} (d\xi) \int_0^\infty
\mathbf{1}_B (r\xi) \nu_\xi(dr), \quad B \in \mathcal{B} (\R^d
\setminus \{ 0 \}).$$ Then $\nu$ is a L\'evy measure, since
\begin{align*}
\int_{S}\wt{\ld}(d\xi)&\int_0^\infty (r^2 \wedge 1)\nu_\xi(dr)\\
& \leq  \int_S \wt{\ld} (d\xi) \int_0^1 r^2 \, \nu_\xi(dr) +
\int_S \wt{\ld} (d\xi) \int_1^\infty \nu_\xi(dr) \\
& =  \int_S \wt{\ld} (d\xi) \int_1^\infty \al^{-1}r^{-2/\al} r^{-1} \wt{Q}_\xi (dr)
+ \int_S \wt{\ld} (d\xi)\int_0^1 \al^{-1} r^{-1} \wt{Q}_\xi (dr),
\end{align*}
which is finite by \eqref{eq-Q2}. If $\mu$ is any infinitely
divisible distribution with L\'evy measure $\nu$, then part (i) of
the proof shows that $\mE_\alpha (\mu)$ has the given L\'evy measure
$\wt{\nu}$, and from the transformation of the generating
triplet in Proposition~\ref{prop-1} we see that $\mu_0 \in
I(\R^d)$ can be chosen such that $\mE_\alpha (\mu_0) = \wt{\mu}$.
 \qed

\vskip 10mm

%%%%%%%%%%%%%%%%%%%%% Section 4 %%%%%%%%%%%%%%%%%%%%%%%%%%%%
\section{The composition of $\Phi$ with $\mE_\alpha$}
In this section we study the composition $\Phi\circ \mE_\alpha$. We
start with the following proposition.

\begin{prop} \label{prop-2}
Let $\alpha > 0$, $m\in \{1,2,\ldots \}$ and $\mu \in I(\R^d)$. Then
$\mu \in I_{\log^m}(\R^d)$ if and only if $\mE_\alpha (\mu) \in
I_{\log^m}(\R^d)$.
\end{prop}

\noindent {\it Proof.} Let $\nu$ and $\wt{\nu}$ denote the L\'evy
measures of $\mu$ and $\mE_\alpha (\mu)$, respectively. By
\eqref{eq-tripel-nu}, we conclude that
$$\int_{\R^d} \varphi(x) \, \wt{\nu} (dx) = \int_{\R^d} \nu(dx)\int_0^\infty
\varphi(ux) \alpha u^{\alpha -1} e^{-u^\alpha}\, du\, $$ for every
measurable nonnegative function $\varphi:\R^d \to [0,\infty]$. In
particular, we have
\begin{align*}
\int_{|x|>1} (\log |x|)^m \, \wt{\nu} (dx) & =  \int_{\R^d}\nu(dx)
\int_{1/|x|}^\infty (\log (u|x|))^m \, \alpha
u^{\alpha-1} e^{-u^\alpha} \, du \, \\
& =  \int_{\R^d} \nu(dx) \sum_{n=0}^m {m \choose n} (\log |x|)^{m-n}
\int_{1/|x|}^\infty (\log u)^n  \alpha u^{\alpha-1} e^{-u^\alpha} \,
du \\
& =:  \int_{\R^d} h(x) \nu(dx), \quad \mbox{say.}
\end{align*}
Then it is easy to see that $h(x) = o(|x|^2)$ as $|x|\downarrow 0$
and that $\lim_{|x|\to\infty} h(x)/(\log|x|)^m = \int_0^\infty
\alpha u^{\alpha -1} e^{-u^{\alpha}} \, du = 1$. Hence,
$\int_{|x|>1} (\log|x|)^m \, \wt{\nu}(dx) < \infty$ if and only if
$\int_{|x|>1} (\log|x|)^m \linebreak \nu(dx) < \infty$, giving the
claim. \qed
%%%This was the old proof%%%
%In particular, we have
%\begin{align*}
%\int_{|x|>1} \log |x| \, \wt{\nu} (dx) & =  \int_{\R^d}\nu(dx)
%\int_{1/|x|}^\infty \log (u|x|) \, \alpha
%u^{\alpha-1} e^{-u^\alpha} \, du \, \\
%& =  \int_{\R^d} \left( (\log |x|) \, e^{-2^\alpha/|x|^\alpha} +
%\int_{1/|x|} (\log u) \, \alpha u^{\alpha-1} e^{-u^\alpha} \, du
%\right) \nu(dx) \\
%& =:  \int_{\R^d} h(x) \nu(dx), \quad \mbox{say.}
%\end{align*}
%Then it is easy to see that $h(x) = o(|x|^2)$ as $|x|\downarrow 0$
%and $\lim_{|x|\to\infty} h(x)/\log|x| = 1$. Hence, $\int_{|x|>1}
%\log|x| \, \wt{\nu}(dx) < \infty$ if and only if $\int_{|x|>1}
%\log|x| \, \nu(dx) < \infty$, giving the claim. \qed

\begin{thm}\label{theorem:th3}
Let $\alpha > 0$ and
$$
 n_\alpha (x)
= \int_x^\infty u^{-1} e^{-u^\alpha} \, du, \q x>0.
$$
Let
$x=n_\alpha^*(t), t>0,$ be its inverse function, and define
the mapping $\mathcal{N}_\alpha : I_{\log}(\R^d) \to I(\R^d)$ by
$$
\mathcal{N}_\alpha (\mu) = \law \left( \int_0^\infty n_\alpha^* (t)
\, dX_t^{(\mu)} \right), \quad \mu\in I_{\log}(\rd).
$$
It then holds
\begin{equation}\label{41a}
\Phi \circ \mE_\alpha = \mE_\alpha \circ \Phi = \mathcal{N}_\alpha,
\end{equation}
including the equality of the domains.
In particular, we have
\begin{align}\label{41}
\Phi \circ\mE _2=\mE_2 \circ \Phi = \mM.
\end{align}
\end{thm}

\begin{rem}
A more general mapping than $\mathcal N_{\al}$-mapping is already
defined in \cite{MN08} and it is shown that $\frak D (\mathcal
N_{\al})=I_{\log}(\rd)$ in Theorem 2.4 of \cite{MN08}.
\end{rem}

\noindent {\it Proof of Theorem~\ref{theorem:th3}.} We remark that
the equation $\Phi\circ \mE_\alpha = \mE_\alpha \circ \Phi$
including the equality of domains can be concluded from
Proposition~\ref{prop-2} and  the general theory of Upsilon
transformations in \cite{BNRT08}, which could also be used to show that
$\Phi\circ \mE_\alpha$ and $\mathcal{N}_\alpha$ transform the L\'evy
measure of the underlying $\mu$ in the same way. To obtain the
transformation of the generating triplet, however, we
give the following proof, which does not refer to the general theory
of Upsilon transformations.

It follows from Proposition \ref{prop-2} that both $\Phi\circ
\mE_\alpha$ as well as $\mE_\alpha \circ \Phi$ are well defined on
$I_{\log} (\R^d)$ and that they have the same domain. Note that
\begin{align*}
C_{\mE_\alpha(\mu)}(z)&=\int_0^1C_{\mu}\left((\log
t^{-1})^{1/\alpha}z\right)dt\end{align*} and
\begin{align*}
C_{\Ph(\mu)}(z)&=\int_0^{\infty}C_{\mu}\left(e^{-t}z\right)dt .
\end{align*}
Then, if we are allowed to exchange the order of the integrals by Fubini's theorem, we have
\begin{align*}
C_{(\mE_\alpha\circ\Ph)(\mu)}(z)
&=\int_0^1dt\int_0^{\infty}C_{\mu}\left((\log
t^{-1})^{1/\alpha}e^{-s}z\right)ds\\
&=\int_0^1dt\int_0^1C_{\mu}\left((\log t^{-1})^{1/\alpha}uz\right)u^{-1}du\\
&=\int_0^1u^{-1}du\int_0^1C_{\mu}\left((\log t^{-1})^{1/\alpha}uz\right)dt\\
&=\int_0^1u^{-1}du\int_0^{\infty}C_{\mu}(vuz)\alpha v^{\alpha-1}
e^{-v^\alpha}dv\\
&=\int_0^{\infty}\alpha v^{\alpha-1} e^{-v^\alpha}dv\int_0^1C_{\mu}(vuz)u^{-1}du\\
&=\int_0^{\infty}\alpha v^{\alpha-1} e^{-v^\alpha}dv\int_0^vC_{\mu}(sz)s^{-1}ds\\
&=\int_0^{\infty}C_{\mu}(sz)s^{-1}ds\int_s^{\infty}\alpha v^{\alpha-1}e^{-v^\alpha}dv\\
&=\int_0^{\infty}C_{\mu}(sz)s^{-1}e^{-s^\alpha}ds\\
&=\int_0^{\infty}C_{\mu}(n_\alpha^{\ast}(t)z)dt,
\end{align*}
and the same calculation can be carried out for $C_{(\Ph\circ
\mE_\alpha)(\mu)}(z)=\int_0^{\infty}C_{\mu}(n_\alpha^{\ast}(t)z)dt$.

 In order to assure the exchange of the order of the
integrations by Fubini's theorem, it is enough to show that
\begin{equation}
\int_0^1u^{-1}du\int_0^{\infty}|C_{\mu}(vuz)| v^{\alpha-1}
e^{-v^\alpha}dv<\infty.
\end{equation}
We have
$$
|C_{\mu}(z)|\leq {2}^{-1}({\rm tr}A)|z|^2 +|\gm||z|+\int_{\rd}|g(z,x)|\nu(dx),
$$
where
$$
g(z,x)=e^{i\langle z,x \rangle}-1-i\langle z,x \rangle (1+|x|^2)^{-1}.
$$
Hence
\begin{align*}
|C_{\mu}(uvz)|&\leq {2}^{-1}({\rm tr}A)u^2v^2|z|^2 +|\gm||u||v||z|
+\int_{\rd}|g(z,uvx)|\nu(dx)\\
&+\int_{\rd}|g(uvz,x)-g(z,uvx)|\nu(dx)=:I_1+I_2+I_3+I_4,
\end{align*}
say. The finiteness of $\int_0^1 u^{-1} du \int_0^\infty (I_1
+I_2)v^{\alpha -1} e^{-v^\alpha} dv$ is trivial. Noting that
$|g(z,x)|\leq c_z|x|^2(1+|x|^2)^{-1}$ with a positive constant $c_z$
depending on $z$, we have
\begin{align*}
\int_0^1 &u^{-1} du \int_{0}^{\infty} I_3 v^{\al -1}e^{-v^{\al}} dv \\
& \leq c_z\int_{\rd}\nu(dx)\int_{0}^{1} u^{-1} du
\int_{0}^{\infty} \frac{(uv|x|)^2}{1+(uv|x|)^2}v^{\al-1}e^{-v^{\al}} dv\\
&=c_z\left(\int_{|x|\leq1}\nu(dx)+\int_{|x|>1}\nu(dx)\right)
\int_{0}^{1}u^{-1} du \int_{0}^{\infty} \frac{(uv|x|)^2}{1+(uv|x|)^2}v^{\al-1}e^{-v^{\al}} dv\\
&=:I_{31}+I_{32},
\end{align*}
say, and
\begin{align*}
I_{31} &\leq c_z\int_{|x|\leq1} |x|^2 \nu(dx) \int_{0}^{1} u\, du
\int_{0}^{\infty} v^{\al +1} e^{-v^{\al}} dv <\infty.
\end{align*}
As to $I_{32}$, we have
\begin{align*}
I_{32} & =  c_z\int_{|x|>1}\nu(dx)
\left ( \int_{0}^{1/|x|^2} + \int_{1/|x|^2}^1\right )u^{-1} du
\int_{0}^{\infty} \frac{(uv|x|)^2}{1+(uv|x|)^2}v^{\al-1}e^{-v^{\al}} dv\\
& =: I_{321}+I_{322},
\end{align*}
say, and
\begin{align*}
I_{321} & \le c_z\int_{|x|>1}\nu(dx)\int _0^{1/|x|^2} du \sek v^{\al +1}e^{-v^{\al}}dv <\infty,\\
I_{322} & \le c_z\int_{|x|>1}\nu(dx) \int _{1/|x|^2}^1 u^{-1}ds \sek v^{\al-1}e^{-v^{\al}}dv\\
& = 2 c_z \int_{|x|>1}\log |x|\nu(dx) \sek
v^{\al-1}e^{-v^{\al}}dv<\infty,
\end{align*}
because $\mu\in I_{\log}(\rd)$.
As to  $I_4$, note that for $a\in\mathbb R$,
\begin{align*}
|g(az,x)- &g(z,ax)|=
\frac{|\langle az,x \rangle ||x|^2|1-a^2|}{(1+|x|^2)(1+|ax|^2)}\\
&\leq \frac{|z||x|^3(|a|+|a|^3)}{(1+|x|^2)(1+|ax|^2)}\\
&\leq
\frac{|z||x|^2(1+|a|^2)}{2(1+|x|^2)},\quad (\text{since}\,\, |b|(1+b^2)^{-1}\le 2^{-1}).\\
& \leq
\frac12 |z|(1+|a|^2).
\end{align*}
Then
\begin{align*}
\int_{0}^{1}& u^{-1} du  \int_{0}^{\infty} I_4 v^{\al-1}e^{-v^{\al}}dv\\
& =  \left (\int_{|x|\le 1}\nu(dx) + \int_{|x|>1}\nu(dx)\right )
 \int_{0}^{1} u^{-1} du \int_{0}^{\infty} |g(uvz,x) - g(z,uvx)| v^{\al-1}e^{-v^{\al}}dv\\
& =:I_{41}+I_{42},
\end{align*}
say.
We have
\begin{align*}
I_{41} & \le |z| \int_{|x|\le 1}\nu(dx)\int_0^1u^{-1}du\sek uv(1+(uv)^2)|x|^3v^{\al-1}e^{-v^{\al}}dv\\
& \le |z| \int_{|x|\le 1} |x|^3 \nu(dx)\sek
(1+v^2)v^{\al}e^{-v^{\al}}dv < \infty.
\end{align*}
Also,
\begin{align*}
I_{42} & \le |z| \int_{|x|> 1}\nu(dx)\left ( \int _0^{1/|x|^3} +
\int_{1/|x|^3}^1\right ) u^{-1}du\sek
\frac {vu(1+(vu)^2)|x|^3}{(1+|x|^2)(1+|vux|^2)}v^{\al-1}e^{-v^{\al}}dv\\
& =: I_{421}+I_{422},
\end{align*}
say.
We have
\begin{align*}
I_{421} & \le |z| \int_{|x|>1}|x|^3 \nu(dx) \int_0^{1/|x|^3}du \sek (1+v^2) v^{\al}e^{-v^{\al}}dv\\
&= |z| \int_{|x|>1}\nu(dx) \sek (1+v^2)v^{\al}e^{-v^{\al}}dv
<\infty,
\end{align*}
and
\begin{align*}
I_{422} & \le \frac{|z|}{2} \int_{|x|>1}\nu(dx)\int_{1/|x|^3}^1 u^{-1}du \sek (1+v^2)v^{\al-1}e^{-v^{\al}}dv\\
& \le \frac{3}{2} |z| \int_{|x|>1}\nu(dx) \log |x| \int_0^\infty
(1+v^2) v^{\alpha-1} e^{-v^\alpha} dv<\infty,
\end{align*}
because $\mu\in I_{\log}(\rd)$.
This completes the proof of (4.3).

By the absolute convergence of the
above integrals, we see that $\int_0^\infty n_\alpha^* (t) \,
dX_t^{(\mu)}$ is indeed definable for every $\mu \in
I_{\log}(\R^d)$ and that
$$C_{(\Phi \circ \mE_\alpha)(\mu)}(z) =
C_{(\mE_\alpha \circ \Phi)(\mu)}(z) = \int_0^\infty C_\mu
(n_\alpha^{\ast}(t) z) \, dt =  C_{\mathcal{N}_\alpha(\mu)}(z),
\quad z\in \R^d,
$$
(see Sato \cite{S07}, Theorem 3.5), and we must have $\Phi\circ
\mE_\alpha = \mE_\alpha \circ \Phi = \mathcal{N}_\alpha$. Since
$\mathcal{N}_2 = \mathcal{M}$, this shows in particular \eqref{41}.
%that $\mathcal{M} = \mE_2 \circ \Phi = \Phi \circ \mE_2$.
\qed

\vskip 3mm
An immediate consequence of Theorem 4.2 is the following.

\begin{thm} \label{theorem-th4}
Let $\alpha >0$. Then
$$
\Phi (E_\alpha(\R^d)\cap I_{\log}(\R^d))
= \mE_\alpha (L(\R^d)) = \mathcal{N}_\alpha (I_{\log}(\R^d)).
$$
\end{thm}

 We
conclude this section with an application of the relation (4.1)
to characterize the limit of certain subclasses obtained by
the iteration of the mapping $\mN$. We need some lemmas. In
the following, $\mN^m$ is defined recursively as $\mN^{m+1} = \mN^m
\circ \mN$.

\begin{lem}\label{lem42}
Let $\al>0.$
For $m=1,2,\ldots$, we have
\begin{equation*}
\frak D(\mN^m)=I_{\log^m}(\rd) \quad \mbox{and}\quad \mN^m = \Phi^m
\circ \mE_{\al}^m = \mE_{\al}^m \circ \Phi^m.
\end{equation*}
\end{lem}

\begin{proof}
%For the proof, we follow that of Lemma 3.8 in [MS08].
By Proposition~\ref{prop-2}, we have  $\mu\in I_{\log^m}(\rd)$ if
and only if $\mathcal E_{\al}(\mu)\in I_{\log^m}(\rd)$.
%We first show that (i) $\mu\in I_{\log^m}(\rd)$ if and only if $\mathcal E_2(\mu)\in I_{\log^m}(\rd)$.
%\lq\lq If part'' is proved as follows.
%If $\int_{|y|>1}(\log |y|)^m\nu_{\mu}(dy)<\infty$, then
%\begin{align*}
%\int_{|x|>1}&(\log |x|)^m\nu_{\mathcal E_2(\mu)}(dx)=\int_{|x|>1}(\log |x|)^m\int_0^{\infty}
%\nu_{\mu}(s^{-1}dx)se^{-s^2}ds\\
%&=\int_{\rd}\nu_{\mu}(dy)\int_0^{\infty}(\log |sy|)^mse^{-s^2}1_{\{|sy|>1\}}ds\\
%&=\int_{|y|>0}\nu_{\mu}(dy)\int_{1/|y|}^{\infty}(\log |sy|)^mse^{-s^2}ds\\
%&=\int_{|y|>0}\nu_{\mu}(dy)\int_{1/|y|}^{\infty}\sum_{n=0}^m
%{m\choose n} (\log s)^n(\log |y|)^{m-n}se^{-s^2}ds\\
%&\le\sum_{n=0}^m {m \choose n} \int_{|y|>1}(\log |y|)^{m-n}\nu_{\mu}(dy)
%\int_{1/|y|}^{\infty}(\log s)^nse^{-s^2}ds+C\\
%&\le\sum_{n=0}^m {m\choose n}\int_{|y|>1}(\log |y|)^{m-n}\nu_{\mu}(dy)
%\int_{1/|y|}^1(\log s)^nse^{-s}ds\\
%&\q +\sum_{n=0}^m {m\choose n}\int_{|y|>1}(\log |y|)^{m-n}\nu_{\mu}(dy)
%\int_1^{\infty}(\log s)^nse^{-s^2}ds+C<\infty,
%\end{align*}
%where $C>0$.
%\lq\lq Only if part'' can be shown by the same as for that of Lemma 3.8 in [MS08].
%So we omit it here.
As shown in the proof of Lemma~3.8 in~\cite{MS08}, we also have that
$\mu\in I_{\log^{m+1}}(\rd)$ if and only if $\mu\in I_{\log}(\rd)$
and $\Phi(\mu)\in I_{\log^m}(\rd)$, and thus $\frak D (\Phi ^m) =
I_{\log ^m}(\rd)$. Since $\mathcal N _{\al}= \Phi \circ\mE _{\al} =
\mE_{\al}\circ \Phi$, we conclude that \begin{equation}
\label{eq-iii} \mbox{$\mu\in I_{\log^{m+1}}(\rd)$ if and only if
$\mu\in I_{\log}(\rd)$ and $\mN(\mu)\in
I_{\log^m}(\rd)$.}\end{equation} Now we prove $\frak
D(\mN^m)=I_{\log^m}(\rd)$ inductively. For $m=1$ this is known, so
assume that $\frak D(\mN^m)=I_{\log^m}(\rd)$ for some $m\ge 1$. If
$\mu\in\frak D(\mN^{m+1})$, then $\mN^{m+1}(\mu)=\mN^m(\mN(\mu))$ is
well-defined. Thus, $\mN(\mu)\in\frak D(\mN^m)=I_{\log^m}(\rd)$ by
assumption, so that $\mu\in I_{\log^{m+1}}(\rd)$ by \eqref{eq-iii}.
Conversely, if $\mu\in I_{\log^{m+1}}(\rd)$, then $\mu\in
I_{\log}(\rd)$ and $\mN(\mu)\in I_{\log^m}(\rd)$ by \eqref{eq-iii},
so that $\mN^m(\mN(\mu))$ is well-defined by assumption. This shows
$\frak D(\mN^{m+1})=I_{\log^{m+1}}(\rd)$. That $\mN^m = \Phi^m \circ
\mE_{\al}^m = \mE_{\al}^m \circ \Phi^m$ for every $m$ then follows
easily from \eqref{41a}, Proposition~\ref{prop-2} and $\frak D (\Phi
^m) = I_{\log^m}(\rd)$.
\end{proof}

Let $S(\rd)$ be the class of all stable distributions on $\rd$, and
for $m=0,1,\ldots$ denote $L_m(\rd) =
\Phi^{m+1}(I_{\log^{m+1}}(\rd))$, $L_{\infty}(\rd)=\cap
_{m=0}^{\infty}L_m(\rd), N_{\al, m}(\rd) =
\mN^{m+1}(I_{\log^{m+1}}(\rd))$ and $N_{\al, \infty}(\rd)=\cap
_{m=0}^{\infty}N_{\al,m}(\rd)$. It is known (cf. Sato, \cite{S80}) that
$L_{\infty}(\rd) = \overline{S(\rd)}$, where the closure is taken
under weak convergence and convolution.
%On the other hand,
%Aoyama~\cite{A08} has recently proved that $M_{\infty}(\rd)\cap
%I_{\rm sym}(\rd) \supset S(\rd)\cap I_{\rm sym}(\rd)$. The same
%proof works for showing that \begin{equation} \label{eq-supset}
%M_{\infty}(\rd)\supset S(\rd).\end{equation}
In order to show that also $N_{\al,\infty} (\bR^d) = \overline{S(\rd)}$,
we need two further lemmas.

\begin{lem} \label{lem-45a}
For $\alpha >0$, $\mE_\alpha$ maps $S(\bR^d)$ bijectively onto
$S(\bR^d)$, namely
$$\mE_\alpha ( S(\rd)) = S(\rd).$$
\end{lem}
This is an immediate consequence of Proposition~\ref{prop-1}~(ii).

\begin{lem}\label{lem46}
Let $\al >0$.
For $m=0,1,\dots$, $N_{\al,m}(\rd)$ is closed under convolution and weak
convergence, and \begin{equation} \label{eq-4.5} S(\rd) \subset
N_{\al,m}(\rd) = \mE_{\al}^{m+1} ( L_m (\rd)) \subset L_m (\rd).\end{equation}
\end{lem}

\begin{proof}
By Lemma $\ref{lem42}$,
$$N_{\al,m}(\rd) = \mN^{m+1} ( I_{\log^{m+1}} (\rd)) = (\mE_{\al}^{m+1} \circ
\Phi^{m+1}) (I_{\log^{m+1}} (\rd))=\mathcal E_{\al}^{m+1}(L_m(\rd)),
$$
hence $S(\rd) \subset N_{\al,m}(\bR^d)$ by Lemma~\ref{lem-45a} and the
fact that $S(\rd) \subset L_m(\rd)$. Further,
$$N_{\al,m}(\rd) = (\Phi^{m+1} \circ \mE_{\al}^{m+1}) (I_{\log^{m+1}} (\rd))
\subset \Phi^{m+1} ( I_{\log^{m+1}} (\rd)) = L_m(\rd).$$ Next
observe that $\mE_{\al}$ and hence $\mE_{\al}^{m+1}$ clearly respect
convolution. Since $L_m(\rd)$ is closed under convolution and weak
convergence (see the proof of Theorem~D in~\cite{BNMS06}), it
follows from \eqref{eq-4.5} and Proposition~\ref{prop-1}~(iv) that
$N_{\al,m}(\rd)$ is closed under convolution and weak convergence, too.
\end{proof}

We can now characterize $N_{\al,\infty} (\rd)$ as the closure of $S(\rd)$
under convolution and weak convergence:

\begin{thm}
Let $\al >0$.
 It holds
$$
 L_{\infty}(\rd) = N_{\al,\infty}(\rd) =\overline {S(\rd)} .
$$
\end{thm}

\begin{proof}
By %\eqref{eq-supset} and
\eqref{eq-4.5} we have
$$
\overline {S(\rd)} = L_{\infty}(\rd) \supset N_{\al,\infty}(\rd) \supset S(\rd).
$$
But since each $N_{\al,m}(\rd)$ is closed under convolution and weak
convergence, so must be the intersection $N_{\al,\infty}(\rd) =
\bigcap_{m=0}^\infty N_{\al,m}(\rd)$, and the assertion follows.
\end{proof}

\vskip 10mm
%%%%%%%%%%%%%%%%%%%%%%%%%%%% Section 5%%%%%%%%%%%%%%%%%%%%%%%%%
\section{Characterization of subclasses of $E _{\al}(\R ^1)$  by stochastic integrals
with respect to some compound Poisson processes} \label{S5}

For any L\'evy process $Y=\{Y_t\}_{t\geq 0}$, denote by
$\mathbf{L}_{(0,\infty)} (Y)$ the class of locally $Y$-integrable
functions on $(0,\infty)$ (cf. Sato~\cite{S07}, Definition~2.3), and
let
\begin{align*}
{\rm Dom} (Y) & = \left\{ h\in \mathbf{L}_{(0,\infty)} (Y)\,\, :\,\,
 \sek
h(t)dY_t\,\,\text{is definable}\right\}\, ,\\
{\rm Dom}^{\downarrow}(Y) & = \{ h\in {\rm Dom} (Y) :
h\,\, \text{is a left-continuous and decreasing function}\\
&  \quad \quad \quad \quad \quad \quad \quad \text{such that}\,\,
\lim_{t\to\infty} h(t)=0\}.
\end{align*}
Here, following Sato~\cite{S07}, Definition~3.1, by saying that the
(improper stochastic integral) $\sek h(t) dY_t$ is definable we mean
that $ \int_p^q h(t) dY_t$ converges in probability as $p\downarrow
0$, $q\to \infty$, with the limit random variable being denoted by
$\sek h(t) dY_t$.

 The property of $h$ belonging to ${\rm Dom}
(Y)$ can be characterized in terms of the generating triplet
$(A_Y,\nu_Y,\gamma_Y)$ of $Y$ and simply properties of $h$, cf.
Sato~\cite{S07}, Theorems~2.6, 3.5 and 3.10. In particular, if
$A_Y=0$, then $h\in {\rm Dom} (Y)$ if and only if $h$ is measurable,
\begin{equation} \label{eq-definable1}
\int_0^\infty ds \int_{\mathbb{R}} (|h(s) x|^2\wedge 1) \, \nu_Y(dx)
< \infty, \end{equation}
\begin{equation} \label{eq-definable2}
 \int_p^q \left| h(s) \gamma_Y +
\int_{\mathbb{R}} h(s) x \left( \frac{1}{1+|h(s) x|^2} -
\frac{1}{1+|x|^2}\right) \nu_Y(dx)\right| ds < \infty
\end{equation} for all $0 < p < q < \infty$ and
\begin{equation} \label{eq-definable3} \lim_{p\downarrow 0, q \to
\infty} \int_p^q \left( h(s) \gamma_Y + \int_{\mathbb{R}} h(s) x
\left( \frac{1}{1+|h(s) x|^2} - \frac{1}{1+|x|^2}\right)
\nu_Y(dx)\right) ds\quad \mbox{exists in $\bR$.}
\end{equation}
In this case, $\int_0^\infty h(t) \, dY_t$ is infinitely divisible
without Gaussian part and its L\'evy measure $\nu_{Y,h}$ is given by
\begin{equation} \label{eq-definable4}
\nu_{Y,h} (B) = \int_0^\infty ds \int_{\bR} \mathbf{1}_B (h(s) x) \,
\nu_Y(dx), \quad B\in \mathcal{B} (\bR \setminus \{0\}).
\end{equation}
 If $\nu_Y$ is symmetric and $\gamma_Y=0$, then
\eqref{eq-definable2} and \eqref{eq-definable3} are automatically
satisfied, so that $h\in {\rm Dom} (Y)$ if and only if
\eqref{eq-definable1} is satisfied, in which case $\gm_Y$  in the
generating triplet of $\int_0^\infty h(t) \,dY_t$ is 0.

Recall the definitions of $E_\alpha^0(\bR^1)$ and $E_\alpha^{0, \rm
sym}(\bR^1)$ from \eqref{def-0} and \eqref{def-sym}. The next
theorem characterizes $E_\alpha^{0,\rm sym} (\bR^1)$ as the class of
distributions which arise as improper stochastic integrals over
$(0,\infty)$ with respect to some fixed symmetric compound Poisson
process.

\begin{thm} \label{theorem:th7}
Let $\alpha >0$ and denote by $Y^{(\alpha)} = \{
Y_t^{(\alpha)}\}_{t\geq 0}$ a compound Poisson process on $\bR$ with
L\'evy measure $\nu_{Y^{(\alpha)}}(dx) =|x|^{\alpha-1}
e^{-|x|^\alpha}dx$ $($without drift$)$. Then
\begin{align}
E _\alpha^{0,\rm sym} (\R^1) &  = \left\{\law \left ( \sek
h(t)dY_t^{(\al)} \right ) \,\,:\,\, h\in
{\rm Dom} (Y^{(\alpha)}) \right\} \label{eq-2}\\
& = \left\{\law \left ( \sek h(t)dY_t^{(\alpha)} \right )
 \,\,:\,\, h\in {\rm Dom}^\downarrow (Y^{(\al)}) \right \} \label{eq-1}.
\end{align}
\end{thm}

\vskip 3mm
\begin{proof}
Let $\mu\in E_\alpha^{0,\rm sym}(\R^1)$. By definition, the L\'evy
measure $\nu$ of $\mu$ has the polar decomposition $(\lambda, \nu_\xi)$
given by
\begin{align} \label{eq-polar1}
\nu_{\xi}(dr)=r^{\alpha-1}g_{\xi}(r^\alpha)dr, \,\, r>0, \,\,
\xi\in\{-1,1\},
\end{align}
and
$$\ld (d\xi) = (\delta _{\{-1\}} + \delta_{\{1\}})(d\xi),$$
where $g_1=g_{-1}$  are completely monotone and $\delta_x$ denotes
Dirac measure at $x$ (If $\mu=\delta_0$ we define $g_\xi=0$ and
shall also call $(\lambda,\nu_\xi)$ a polar decomposition, even if
$\nu_\xi$ is not strictly positive here). In the following, we drop
the subscript $\xi$ of $g_{\xi}$ and $\nu_{\xi}$. Since $g$ is completely monotone,
there exists a Borel measure $Q$ on $[0,\infty)$ such that
$g(y)=\int_{[0,\infty)} e^{-yt}Q(dt)$. By \eqref{eq-Q2}, in order
for $\nu$ to satisfy $ \sek (x^2 \wedge 1)\nu(dx) <\infty$, it is
necessary and sufficient that
\begin{equation} \label{eq-Q}
Q(\{ 0 \}) = 0 , \quad \int_0^1 t^{-1} \, Q(dt) < \infty \quad
\mbox{and}\quad \int_1^\infty t^{-1-2/\alpha} \, Q(dt) < \infty.
\end{equation}
Observe that under this condition, we have for each $r>0$,
\begin{align*} \label{eq-tail1}
\nu([r,\infty))&=\int_r^{\infty}y^{\alpha -1}g(y^\alpha)dy=
\int_0^{\infty}(\alpha t)^{-1}Q(dt) \int_r^{\infty}\alpha t
y^{\alpha -1} e^{-y^\alpha t}dy \\
& =\int_0^{\infty}(\alpha t)^{-1}e^{-r^\alpha t}Q(dt).\nonumber
\end{align*}

Next, observe that $Y^{(\alpha)}$ is symmetric without Gaussian
part, so that by  \eqref{eq-definable1}  a measurable function $h$
is in ${\rm Dom} (Y^{(\alpha)})$ if and only if
\begin{equation} \label{eq-conv}
\int_0^\infty ds \int_\R  \left( |h(s) x|^2 \wedge 1\right)
|x|^{\alpha -1} \, e^{-|x|^\alpha} \, dx \,  < \infty,
\end{equation}
in which case $\int_0^\infty h(t) \, dY_t^{(\alpha)}$ is infinitely
divisible with the generating triplet \linebreak $(A_{Y,h}=0,
\nu_{Y,h}, \gamma_{Y,h}=0)$ and the L\'evy measure  ${\nu}_{Y,h}$ is
symmetric and by \eqref{eq-definable4} satisfies
\begin{align} \label{eq-finiteness}
{\nu}_{Y,h}([r,\infty))&=\int_0^{\infty}ds\int_{r/|h(s)|}^{\infty}x^{\alpha
-1} e^{-x^\alpha}dx =\alpha^{-1}
\int_0^{\infty}e^{-r^\alpha/|h(s)|^\alpha}ds
\end{align}
for every $r>0$.
 Hence, in order to prove \eqref{eq-2} and
\eqref{eq-1}, it is enough
to prove the following:\\
(a) For each Borel measure $Q$ on $[0,\infty)$ satisfying
\eqref{eq-Q} there exists a function $h\in {\rm
Dom}^{\downarrow}(Y^{(\alpha)})$ such that
\begin{equation} \label{eq-equality}
\int_0^\infty t^{-1} \, e^{-r^\alpha t} \, Q(dt) = \int_0^\infty e^{-r^\alpha /
|h(s)|^\alpha} \, ds\quad \forall\; r > 0.
\end{equation}
(b) For each $h \in {\rm Dom} (Y^{(\al)})$ there exists a Borel
measure $Q$ on $[0,\infty)$ satisfying \eqref{eq-Q} such that
\eqref{eq-equality} holds.

To show (a), let $Q$ satisfy \eqref{eq-Q}, and denote
$$F(x) := \int_{(0,x]} t^{-1} \, Q(dt), \quad x \in [0,\infty),$$
and by
$$F^\leftarrow (t) = \inf \{ y \geq 0 : F(y) \geq t\}, \quad t\in
[0,\infty),$$ its left-continuous inverse, with the usual
convention  $\inf \emptyset = +\infty$. Now define
\begin{equation*} \label{eq-h}
h=h_Q : (0,\infty) \to [0,\infty), \quad t\mapsto
\left(F^\leftarrow (t)\right)^{-1/\alpha}.
\end{equation*}
Then $h$ is left-continuous, decreasing, and satisfies
$\lim_{t\to\infty} h(t) = 0$. Denote Lebesgue measure on
$(0,\infty)$ by $m_1$, and consider the function
\begin{equation} \label{eq-T}
T:(0,\infty) \to (0,\infty], \quad s \mapsto h(s)^{-\alpha} =
F^\leftarrow (s).
\end{equation}
Then $(T(m_1))_{|(0,\infty)}$, the image measure of $m_1$ under the
mapping $T$, when restricted to $(0,\infty)$, satisfies
\begin{equation} \label{eq-T-Q}
(T(m_1))_{|(0,\infty)} (dt) = t^{-1} \, Q_{|(0,\infty)}(dt).
\end{equation}
Hence it follows that for every $r>0$,
$$
\int_{(0,\infty)} e^{-r^\alpha/h(s)^\alpha} \, m_1(ds) =
\int_{(0,\infty) \cap \{ s : T(s) \not= \infty \} } e^{-r^\alpha
T(s)} \, m_1(ds) = \int_{(0,\infty)} e^{-r^\alpha t} \, (T( m_1))
(dt),
$$
yielding \eqref{eq-equality}. To show \eqref{eq-conv}, namely that
$h \in {\rm Dom}(Y^{(\alpha)})$, observe that
\begin{align*}
\int_0^\infty ds&\int_\R \left( |h(s)x|^2 \wedge 1\right)
|x|^{\alpha -1} e^{-|x|^\alpha} dx \\
& =  2 \int_0^\infty x^{\alpha + 1} e^{-x^\alpha} dx\int_{ \{ s:
h(s) \leq 1/x\} } h(s)^2 \, ds \, + 2\int_0^\infty ds
\int_{1/h(s)}^\infty x^{\alpha -1} e^{-x^\alpha} \, dx \, \\
& =  2\int_0^\infty x^{\alpha + 1} e^{-x^\alpha}dx \int_{ \{ s:
T(s) \geq x^\alpha\} }T(s)^{-2/\alpha}\, ds \,  + 2 \alpha^{-1}
\int_0^\infty
e^{-T(s)} \, ds \\
& =  2\int_0^\infty x^{\alpha + 1} e^{-x^\alpha}dx \int_{\{t \geq
x^\alpha\}} t^{-1-2/\alpha} Q(dt) \, + 2\alpha^{-1} \int_0^\infty
e^{-t} t^{-1} \, Q(dt)
\end{align*}
by \eqref{eq-T-Q}. The second of these terms is clearly finite by
\eqref{eq-Q}. To estimate the first, observe that
\begin{align*}
\int_0^\infty x^{\alpha + 1}& e^{-x^\al}dx\int_{x^\al}^\infty t^{-1-2/\al}\,Q(dt)\\
& \leq \int_1^\infty x^{\alpha+1} e^{-x^\alpha} dx\int_1^\infty
t^{-1-2/\alpha} \, Q(dt) \,  + \int_0^1 x^{\alpha +1}dx
\int_{1}^\infty t^{-1-2/\alpha}\, Q(dt) \, \\
&\quad \hskip5mm + \int_0^1 x^{\alpha +1} dx\int_{x^\alpha}^1 t^{-1-2/\alpha}\,
Q(dt) ,
\end{align*}
and the first two summands are finite by \eqref{eq-Q}, while the
last summand is equal to
$$\int_0^1 t^{-1-2/\alpha}Q(dt)  \int_0^{t^{1/\alpha}} x^{\alpha+1} \,
dx\, = (\alpha+2)^{-1} \int_0^1 t^{1+2/\alpha} t^{-1-2/\alpha} \,
Q(dt)$$ and hence also finite. This shows \eqref{eq-conv} for $h$
and hence (a).

To show (b), let $h\in {\rm Dom}(Y^{(\al)})$ and assume first that
$h$ is nonnegative. Let $T:(0,\infty) \to (0,\infty]$ be defined by
$T(s) = h(s)^{-\alpha}$ as in \eqref{eq-T}, and consider the image
measure $T(m_1)$. Define the measure $Q$ on $[0,\infty)$ by $Q(\{ 0
\})=0$ and equality \eqref{eq-T-Q}. Since $\int_0^\infty h(t) \,
dY_t^{(\alpha)}$ is automatically infinitely divisible with L\'evy
measure ${\nu}_{Y,h}$ given by \eqref{eq-finiteness}, we have as in
the proof of (a) for every $r>0$
$$\int_{(0,\infty)} e^{-r^\alpha t} (\alpha t)^{-1} \, Q(dt) =
\alpha^{-1} \int_0^\infty e^{-r^\alpha/h(s)^\alpha } \, ds =
{\nu}_{Y,h} ( [r,\infty)).$$ In particular, $Q$ must be a Borel
measure and \eqref{eq-equality} holds. Since the left hand side of
this equation converges and the right hand side is known to be the
tail integral of a L\'evy measure, it follows that \eqref{eq-Q} must
hold. Hence we have seen that $\law (\sek h(t) dY_t^{(\alpha)}) \in
E_\alpha^{0,\rm sym}(\mathbb R^1)$ for nonnegative $h\in {\rm Dom}
(Y^{(\alpha)})$. For general $h\in {\rm Dom} (Y^{(\alpha)})$, write
$h = h^+ - h^-$ with $h^+ := h \vee 0$ and $h^- := (-h) \vee 0$.
Then $h^+, h^- \in {\rm Dom} (Y^{(\alpha)})$ by \eqref{eq-conv}, and
Equation~\eqref{eq-definable4} and the discussion following it show
that $\int_0^\infty h(t) dY_t^{(\alpha)}$ has no Gaussian part,
gamma part 0 and satisfies $\nu_{Y,h} = \nu_{Y, h^+} + \nu_{Y,h^-}.$
%from which $\int_0^\infty h(t) dY_t^{(\alpha)} \in E_\alpha^{0,\rm
%sym}$ follows.
The corresponding Borel measure $Q$ is given by $Q = Q^+ +Q^-$,
where $Q^+$ and $Q^-$ are constructed from $h^+$ and $h^-$,
respectively, completing the proof of (b).
\end{proof}

Next, we ask whether every distribution in $E_\alpha^0(\bR^1)$ can
be represented as a stochastic integral with respect to the compound
Poisson process $Z^{(\alpha)}$ having L\'evy measure
$\nu_{Z^{(\al)}}(dx) = x^{\alpha -1} e^{-x^\alpha} {\bf
1}_{(0,\infty)} (x) \, dx$ (without drift) plus some constant. We
shall prove that such a statement is true e.g. for those
distributions in $E_\alpha^0(\bR^1)$ which correspond to L\'evy
processes of bounded variation, but that not every distribution in
$E_\alpha^0(\bR^1)$ can be represented in this way. However, every
distribution in $E_\alpha^0(\bR^1)$ appears as an essential limit of
locally $Z^{(\alpha)}$-integrable functions. Following
Sato~\cite{S07}, Definition~3.2, for a L\'evy process $Y=\{
Y_t\}_{t\geq 0}$ and a locally $Y$-integrable function $h$ over
$(0,\infty)$ we say that {\it the essential improper stochastic
integral on $(0,\infty)$ of $h$ with respect to $Y$ is definable} if
for every $0<p<q<\infty$ there are real constants $\tau_{p,q}$ such
that $\int_p^q h(t) dY_t - \tau_{p,q}$ converges in probability as
$p\downarrow 0$, $q\to\infty$. We write ${\rm Dom}_{\rm es} (Y)$ for
the class of all locally $Y$-integrable functions $h$ on
$(0,\infty)$ for which  the essential improper stochastic integral
with respect to $Y$ is definable, and for each $h \in {\rm Dom}_{\rm
es} (Y)$ we denote the class of distributions arising as possible
limits $\int_p^q h(t) dY_t - \tau_{p,q}$ as $p\downarrow 0$, $q\to
\infty$ by $\Phi_{h,\rm es} (Y)$ (the limit is not unique, since
different sequences $\tau_{p,q}$ may give different limit random
variables). As for ${\rm Dom} (Y)$, the property of belonging to
${\rm Dom}_{\rm es} (Y)$ can be expressed in terms of the
characteristic triplet $(A_Y,\nu_Y,\gamma_Y)$ of $Y$. In particular,
if $A_Y = 0$, then a function $h$ on $(0,\infty)$ is in ${\rm
Dom}_{\rm es} (Y)$ if and only if $h$ is measurable and
\eqref{eq-definable1} and \eqref{eq-definable2} hold, and in that
case $\Phi_{h,\rm es} (Y)$ consists of all infinitely divisible
distributions ${\mu}$ with characteristic triplet $({A}_{Y,h} = 0,
\nu_{Y,h}, \gamma)$, where $\nu_{Y,h}$ is given by
\eqref{eq-definable4} and $\gamma \in \bR$ is arbitrary
(cf.~\cite{S07}, Theorems~3.6 and~3.11).

Denote \begin{align*} E_\alpha^+ (\bR^1) & :=  \{ \mu \in
E_\alpha(\bR^1) : \mu ( ( -\infty,0)) = 0\},\\
E_\alpha^{+,0}(\bR^1) & :=  \{ \mu \in E_\alpha^+ (\bR^1) : \,
\mbox{$\{X_t^{(\mu)}\}$ has zero drift}\},\\
 E^{BV}_\alpha (\bR^1) & :=  \{ \mu \in
E_\alpha(\bR^1) : \, \mbox{$\{X_t^{(\mu)}\}$ is of bounded
variation}\},\\
E_{\alpha}^{BV,0} (\bR^1)  &:=  \{ \mu \in E^{BV}_\alpha(\bR^1): \,
\mbox{$\{X_t^{(\mu)}\}$ has zero drift}\}. \end{align*} We then
have:

\begin{thm} \label{theorem:th8}
Let $\alpha >0$ and denote by $Z^{(\alpha)} = \{
Z_t^{(\alpha)}\}_{t\geq 0}$ a compound Poisson process on $\bR$ with
L\'evy measure $\nu_{Z^{(\al)}}(dx) = x^{\alpha -1} e^{-x^\alpha}
{\bf 1}_{(0,\infty)} (x) \, dx$ $($without
drift$)$. Then it holds:\\
$(i)$ The class of distributions arising as limits of essential
improper stochastic integrals with respect to $Z^{(\alpha)}$ is
$E_\alpha^0(\bR^1):$
\begin{equation} \label{eq-i-1}
E_\alpha^0(\bR^1) = \bigcup_{h \in {\rm Dom}_{\rm es}
(Z^{(\alpha)})} \Phi_{h,\rm es} (Z^{(\alpha)}).
\end{equation}
$(ii$) Distributions in $E_\alpha^{BV,0}(\bR^1)$ and $E_\alpha^{+,0}
(\bR^1)$ can be expressed as improper stochastic integrals over
$(0,\infty)$ with respect to $Z^{(\alpha)}$. More precisely
\begin{align}
\label{eq-ii-0} E_\alpha^{+,0}(\bR^1) & =  \left\{ \law \left( \sek
h(t) dZ_t^{(\alpha)}\right) : h \in {\rm Dom} (Z^{(\alpha)}), h \geq
0 \right\},\\
 \label{eq-ii-1} E_\alpha^{BV,0} (\bR^1) &  =
\left\{ \law \left( \sek h(t) dZ_t^{(\alpha)}\right) : h \in {\rm
Dom} (Z^{(\alpha)}) \; \mbox{\rm
such that }\right.\\
&\left. \quad \quad \q\q \int_0^\infty ds \int_{\bR} (|h(s)x| \wedge
1) \nu_{Z^{(\alpha)}} (dx) < \infty.\right\} \nonumber
\end{align} In particular,
\begin{equation} \label{eq-ii-2}
E_\alpha^+ (\bR^1) = \left\{ \law \left( \sek h(t)
dZ_t^{(\alpha)} + b \right) : h \in {\rm Dom} (Z^{(\alpha)}),\,
h\geq 0, \, b\in [0,\infty)\right\}.
\end{equation}
$(iii)$ Not every distribution in $E_\alpha^0(\bR^1)$ can be
represented as an improper stochastic integral over $(0,\infty)$
with respect to $Z^{(\alpha)}$ plus some constant. It holds
\begin{equation} \label{eq-iii-1} E_\alpha^{BV} (\bR^1) \cup
E_\alpha^{0,\rm sym} (\bR^1) \subsetneqq \left\{ \law \left( \sek
h(t) dZ_t^{(\alpha)} + b\right) : b \in \bR, h \in {\rm
Dom}(Z^{(\alpha)})\right\} \subsetneqq  E_\alpha^0 (\bR^1).
\end{equation}
\end{thm}

\begin{proof}
(i)  Let $h \in {\rm Dom}_{\rm es} (Z^{(\alpha)})$ and $\mu\in
\Phi_{h,\rm es} (Z^{(\alpha)})$ and write $h=h^+-h^-$ with $h^+$ and
$h^-$ being the positive and negative parts of $h$, respectively.
Then $\mu$ is infinitely divisible without Gaussian part and by
\eqref{eq-definable4} its L\'evy measure $\nu_{Z,h}$ satisfies
\begin{align*}
\nu_{Z,h,1} ( [r,\infty)) & :=  \nu_{Z,h} ( [r,\infty)) =
\alpha^{-1} \int_0^\infty e^{-r^\alpha / h^+(s)^\alpha} \, ds ,\\
\nu_{Z,h,-1} ( [r,\infty)) & :=  \nu_{Z,h} ( (-\infty,-r]) =
\alpha^{-1} \int_0^\infty e^{-r^\alpha / h^-(s)^\alpha} \, ds
\end{align*}
for every $r>0$. Define the mappings $T_1, T_{-1}: (0,\infty) \to
(0,\infty]$ by $T_1(s) = (h^+(s))^{-\alpha}$ and $T_{-1}(s) = (h^-
(s))^{-\alpha}$ and the measures $Q_1$ and $Q_{-1}$ on $[0,\infty)$
by
$$Q_\xi (\{0\}) = 0 \quad \mbox{and}\quad (T_\xi (m_1))_{|(0,\infty)} (dt) =
t^{-1} {Q_\xi}_{| (0,\infty)} (dt), \quad  \quad \xi\in \{-1,1\}.$$
Then as in the proof of Theorem~\ref{theorem:th7},
$$\int_{(0,\infty)} e^{-r^{\alpha }t} (\alpha t)^{-1} Q_\xi (dt) =
\nu_{Z,h,\xi} ([r,\infty)), \quad r> 0,\quad \xi \in \{-1,1\},$$ and
$Q_1$ and $Q_{-1}$ satisfy \eqref{eq-Q} and we conclude that
$\nu_{Z,h,\xi}(dr) = r^{\alpha -1} g_\xi(r^\alpha) dr$ for
completely monotone functions $g_1$ and $g_{-1}$, so that
$\Phi_{h,\rm es} (Z^{(\alpha)}) \subset E_\alpha^0(\bR^1)$, giving
the inclusion ``$\supset$'' in equation~\eqref{eq-i-1}.

Now let $\mu \in E_\alpha^0 (\bR^1)$ with L\'evy measure $\nu$, and
define the L\'evy measures $\nu_1$ and $\nu_{-1}$ supported on
$[0,\infty)$ by \begin{equation} \label{eq-hallo1} \nu_1(B) := \nu
(B), \quad \nu_{-1} (B) := \nu (-B), \quad B \in \mathcal{B}
((0,\infty)).\end{equation}
 Then
\begin{equation} \label{eq-hallo2}
\nu_\xi ( [r,\infty)) = \int_0^\infty (\alpha t)^{-1} e^{-r^\alpha
t} Q_\xi (dt), \quad r> 0, \quad \xi \in \{ -1,1\}, \end{equation}
for some Borel measures $Q_1$ and $Q_{-1}$ satisfying \eqref{eq-Q}.
As in the proof of (a) in Theorem~\ref{theorem:th7}, we find
nonnegative and decreasing functions $h_1, h_{-1} : (0,\infty) \to
[0,\infty)$ such that \eqref{eq-conv} (i.e. \eqref{eq-definable1}
with $\nu_{Z^{(\alpha)}}$ in place of $\nu_Y$) and
\eqref{eq-equality} hold. Since $h_1, h_{-1}$ are bounded on compact
subintervals of $(0,\infty)$ and since $Z^{(\alpha)}$ has bounded
variation, it follows that $h_1$ and $h_{-1}$ satisfy also
\eqref{eq-definable2}, so that $h_1, h_{-1} \in {\rm Dom}_{\rm es}
(Z^{(\alpha)})$ and the L\'evy measures of $\widetilde{\mu}_1 \in
\Phi_{h_1,\rm es}(Z^{(\al)})$ and $\widetilde{\mu}_{-1} \in \Phi_{h_{-1}, \rm
es}(Z^{(\al)})$ are given by $\nu_1$ and $\nu_{-1}$, respectively. Now define
the function $h:(0,\infty) \to \bR$ by \begin{equation}
\label{eq-hallo3} h(t) = \left\{
\begin{array}{lll} h_1(t-n) , & t \in
(2n,2n+1],  & n \in \{1,2,\ldots\},\\
-h_{-1}(t-n-1) , & t \in
(2n+1,2n+2],  & n \in \{1,2,\ldots\},\\
h_1 ( t - 2^{-k-1}),\quad { } & t \in (2^{-k}, 2^{-k} + 2^{-k-1}], &
k\in
\{0,1,2,\ldots\},\\
-h_{-1} ( t - 2^{-k}), & t \in (2^{-k}+2^{-k-1}, 2^{-k+1} ],\quad {
} & k\in \{0,1,2,\ldots\}. \end{array} \right.
\end{equation}
Then also $h\in {\rm Dom}_{\rm es} (Z^{(\alpha)})$ and any
$\widetilde{\mu} \in \Phi_{h,\rm es} (Z^{(\alpha)})$ has L\'evy
measure $\nu$, showing the inclusion \lq\lq$\subset$'' in
equation~\eqref{eq-i-1}.

(ii) Let $h \in {\rm Dom} (Z^{(\alpha)})$. Then $\int_0^\infty h(t)
\, dZ_t^{(\alpha)} \in E_\alpha^0 (\bR^1)$ by (a). Further, by
Theorem~3.15 in Sato~\cite{S07}, $\int_0^\infty h(t)
dZ_t^{(\alpha)}$ is the distribution at time 1 of a L\'evy process
of bounded variation if and only if \begin{equation} \label{eq-Q4}
\int_0^\infty ds \int_{\bR} (|h(s)x| \wedge 1) \nu_{Z^{(\alpha)}}
(dx) < \infty, \end{equation}
in which case this L\'evy process will
have zero drift. Since $\law (\int_0^\infty h(t) dZ_t^{(\alpha)})$
has trivially support contained in $[0,\infty)$ if $h\geq 0$, this
gives the inclusion ``$\supset$'' in \eqref{eq-ii-0} and
\eqref{eq-ii-1}.

Now suppose that $\mu \in E_\alpha^{BV,0}(\bR^1)$ with L\'evy
measure $\nu$, define $\nu_1$ and $\nu_{-1}$ by \eqref{eq-hallo1}
and choose Borel measures $Q_1$ and $Q_{-1}$ such that
\eqref{eq-hallo2} holds. Then it can be shown in complete analogy to
the proof leading to \eqref{eq-Q2} that for $\xi \in \{-1,1\}$,
$\nu_\xi$ satisfies $\sek (1 \wedge x) \nu_\xi(dx) < \infty$ if and
only if
\begin{equation} \label{eq-Q3}
Q_\xi (\{0\}) = 0 , \quad \int_0^1 t^{-1} Q_\xi (dt) < \infty \quad
\mbox{and}\quad \int_1^\infty t^{-1 - 1/\alpha} Q_\xi (dt) < \infty.
\end{equation}
For $\xi\in \{-1,1\}$ and $x\in [0,\infty)$ define $F_\xi(x) :=
\int_{(0,x]} t^{-1} Q_\xi(dt)$, $h_\xi =
(F_\xi^\leftarrow)^{-1/\alpha}$ and $T_\xi = (h_\xi)^{-\alpha} =
F_\xi^\leftarrow$. Then it follows in complete analogy to the proof
of (a) of Theorem~\ref{theorem:th7}, using \eqref{eq-Q3}, that
\eqref{eq-equality} and \eqref{eq-Q4} hold for $h_\xi$ and $Q_\xi$.
By Theorem~3.15 in Sato~\cite{S07} this then shows that $h_\xi  \in
{\rm Dom} (Z^{(\alpha)})$ for $\xi\in \{-1,1\}$. Now if $\mu \in
E_\alpha^{+,0} (\bR^1)$, define $h(t) := h_1(t)$, and for general
$\mu \in E_\alpha^{BV,0}$, define $h(t)$ by \eqref{eq-hallo3}. In
each case $h$ satisfies \eqref{eq-Q4}, $h\in {\rm Dom}
(Z^{(\alpha)})$, and $\mu = \law (\int_0^\infty h(t)
dZ_t^{(\alpha)})$, giving the inclusions ``$\subset$'' in
\eqref{eq-ii-0} and \eqref{eq-ii-1}.

(iii)
 Let $\mu \in E_\alpha^{0,\rm sym}(\mathbb R^1)$. By Theorem~\ref{theorem:th7}
there exists $f\in {\rm Dom}^\downarrow (Y^{(\alpha)})$ such that
$\mu = \law ( \sek f(t) dY^{(\alpha)}_t)$. Write $h_1 = h_{-1} := f$
and define the function $h:(0,\infty) \to \bR$ by \eqref{eq-hallo3}.
We claim that $h\in {\rm Dom} (Z^{(\alpha)})$. To see this, observe
that $h$ clearly satisfies \eqref{eq-definable1} with respect to
$\nu_{Z^{(\alpha)}}$ since $f$ has the corresponding property with
respect to $\nu_{Y^{(\alpha)}}$. Next, since ${|h(s) x|}({1+
|h(s) x|^2})^{-1}$ is bounded by $1/2$ and $\nu_{Z^{(\alpha)}}(\bR)$ is
finite, it follows that \begin{equation} \label{eq-hallo5} \int_0^q
\left| \sek \frac{h(s) x}{1+ |h(s) x|^2} x^{\alpha-1}
e^{-x^{\alpha}} dx \right| ds < \infty \quad \forall\; q >
0.\end{equation} But since $Z^{(\alpha)}$ has the generating triplet
$$\left(A_{Z^{(\alpha)}}=0,\,\,\nu_{Z^{(\alpha)}},\,\,
\gamma_{Z^{(\alpha)}} = \sek \frac{x}{1+x^2} x^{\alpha -1}
e^{-x^{\alpha}} dx\right),$$ \eqref{eq-hallo5} shows that
\eqref{eq-definable2} is satisfied for $h$ with respect to
$\nu_{Z^{(\alpha)}}$. Finally, by the definition of $h$, for
$$\gamma_{Z,h,0,q} := \int_0^q \left(\sek \frac{h(s) x}{1 + |h(s) x|^2}
x^{\alpha -1} e^{-x^{\alpha}} \, dx\right) ds, \quad q > 0,$$ we
have $\gamma_{Z,h,0,q} = 0$ for $q=2,4,6,\ldots$, and since
$\lim_{t\to \infty} h(t) = 0$ it follows that $\lim_{q\to\infty}
\gamma_{Z,h,0,q}$ exists and is equal to 0. We conclude that
\eqref{eq-definable3} is satisfied, so that $h \in {\rm Dom}
(Z^{(\alpha)})$. By \eqref{eq-definable4} we clearly have $\law
\left(\sek h(t) dZ_t^{(\alpha)}\right) = \law \left(\sek f(t)
dY_t^{(\alpha)}\right) = \mu.$ Together with \eqref{eq-i-1} and
\eqref{eq-ii-1} and this shows \eqref{eq-iii-1} apart from the fact
that the inclusions are proper.

To show that the first inclusion in \eqref{eq-iii-1} is proper, let
$\mu \in E_\alpha^{0, \rm sym} (\bR^1) \setminus E_\alpha^{BV}
(\bR^1)$. The latter set is nonempty since by (5.8) and
\eqref{eq-Q3} it suffices to find a Borel measure $Q$ on
$[0,\infty)$ such that (5.8) holds but $\int_1^\infty
t^{-1-1/\alpha} Q(dt) = \infty$. As already shown, there exists
$h\in {\rm Dom} (Z^{(\alpha)})$ such that $\mu = \law ( \sek h(t)
dZ_t^{(\alpha)})$. Then $h+ {\bf 1}_{[1,2]} \in {\rm Dom} (Z^{(\alpha)})$,
and $\law (\sek (h(t) + {\bf 1}_{[1,2]} (t) dZ_t^{(\alpha)})$ is clearly
neither symmetric nor of finite variation.

To see that the second inclusion in \eqref{eq-iii-1} is proper, let
$\mu \in E_\alpha^0(\bR^1)$ with L\'evy measure $\nu$ being
supported on $[0,\infty)$ such that $\int_0^1 x \, \nu(dx) =
\infty$. Suppose there are $b\in \bR$ and $h\in {\rm Dom}
(Z^{(\alpha)})$ such that $\mu = \law (\sek h(t) dZ_t^{(\alpha)}
+b)$. Since $\nu$ is supported on $[0,\infty)$, we must have $h\geq
0$ Lebesgue almost surely, so that we can suppose that $h\geq 0$
everywhere. Then we have from \eqref{eq-definable1} and
(5.3) that
$$\sek ds \sek (|h(s) x|^2 \wedge 1) \, \nu_{Z^{(\alpha)}} (dx) <
\infty$$ and
$$\sek ds \sek \frac{h(s) x}{1+h(s) x} \nu_{Z^{(\alpha)}} (dx) <
\infty.$$ Together these two equations imply
$$\sek ds \sek (|h(s) x| \wedge 1) \nu_{Z^{(\alpha)}} (dx) <
\infty,$$ so that $\mu \in E_\alpha^{BV} (\bR^1)$ by
\eqref{eq-ii-1}, contradicting $\int_0^1 x \, \nu(dx) = \infty$.
This completes the proof of \eqref{eq-iii-1}.
\end{proof}

In the following, we shall call a class $F$ of distributions in
$\bR^1$ {\it closed under scaling} if for every $\bR^1$-valued
random variable $X$ such that $\law (X) \in F$ it also holds that
$\law (cX) \in F$ for every $c>0$. If $F$ is a class of infinitely
divisible distributions in $\bR^1$ and satisfies that $\mu \in F$ implies
 $\mu^{s \ast}\in F$ for any $s>0$, where $\mu^{s\ast}$ is the
distribution with characteristic function $(\widehat{\mu}(z))^s$, we
shall call $F$ {\it closed under taking of  powers}. Recall
that a class $F$ of infinitely divisible distributions on $\bR^1$ is
called {\it completely closed in the strong sense} if it is closed
under convolution, weak convergence,  scaling, taking of powers, and
additionally contains $\mu \ast \delta_b$ for any $\mu \in F$ and
$b\in \bR$. We can now characterize $E_\alpha(\bR^1)$ and certain
subclasses as smallest classes which satisfy certain properties.

\begin{thm}  \label{theorem:th9}
Let $\alpha >0$ and
$Y^{(\alpha)}$ and $Z^{(\alpha)}$ be defined as in
Theorems~\ref{theorem:th7} and~\ref{theorem:th8}, respectively. Then it holds$:$\\
$(i)$ The class $E_\alpha(\R^1)$  is the smallest class of infinitely
divisible distributions on $\R^1$ which is completely closed in the
strong sense and contains $\law(Z_1^{(\alpha)})$ and  $\law
(-Z_1^{(\alpha)})$.\\
$(ii)$ The class $E_\alpha^+(\bR^1)$ is the smallest class of
infinitely divisible distributions on $\bR^1$ which is closed under
convolution, weak convergence,  scaling, taking of powers and
contains $\law (Z_1^{(\alpha)})$.\\
$(iii)$ The class $E_\alpha^{\rm sym} (\bR^1):= E_\alpha (\bR^1)\cap
I_{\rm sym} (\bR^1)$ is the smallest class of infinitely divisible
distributions on $\bR^1$ which is closed under convolution, weak
convergence,  scaling, taking of powers and contains $\law
(Y_1^{(\alpha)})$.
\end{thm}

\begin{proof}
{}From the definition it is clear that all the classes under
consideration are closed under convolution, scaling and taking of
powers. The class $E_\alpha(\bR^1)$ is closed under weak convergence
by Proposition~\ref{prop-1}(iv) and Theorem~\ref{theorem:th1}, and
hence so are $E_\alpha^+(\bR^1)$ and $E^{\rm sym}_\alpha (\bR^1)$.
Further, all the given classes contain the specified distributions
which can be seen by taking $h={\bf 1}_{(0,1]}$ (or also $-{\bf
1}_{(0,1]}$ for (i)) in Theorems~\ref{theorem:th7} and
\ref{theorem:th8}, respectively. So it only remains to show that the
given classes are the smallest classes among all classes with the
specified properties.

(iii) Let $F$ be a class of infinitely divisible distributions which
is closed under convolution, weak convergence, scaling, taking of
powers and   which contains $\law(Y_1^{(\alpha)})$. Since $F$ is
closed under taking of powers it contains also $\law
(Y_t^{(\alpha)})$ for all $t>0$, and by the scaling and convolution
property also
$$\law\left( \int_{(0,\infty)} h_n(t) \,
dY_t^{(\alpha)}\right) = \law\left( \sum_{j=1}^{k_n} \beta_{j,n}
(Y_{t_{j,n}}^{(\alpha)} - Y_{t_{j-1,n}}^{(\alpha)})\right)$$ for
every simple function $h_n = \sum_{j=1}^{k_n} \beta_{j,n}
1_{(t_{j-1,n},t_{j,n}]}$ with $\beta_{j,n} \geq 0$ and $0 < t_{1,n}
< \ldots < t_{k_n, n}$. Now let $h$ be an arbitrary function  in
${\rm Dom}^{\downarrow}(Y^{(\alpha)})$ and choose a sequence
$(h_n)_{n\in\mathbb{N}}$ of such simple functions such that $h_n$
converges pointwise from below to $h$. Then $\int_{[a,b]} h_n (t)\,
dY_t^{(\alpha)}$ converges in probability (even almost surely) to
$\int_{[a,b]} h(t)\, dY_t^{(\alpha)}$ for all $0<a<b<\infty$, from
which it follows that $\law \left( \int_{[a,b]} h(t) \,
dY_t^{(\alpha)}\right)\in F$ for all $0<a<b<\infty$ and hence $\law
\left( \int_0^\infty h(t) dY_t^{(\alpha)}\right) \in F$.  Together
with Theorem~\ref{theorem:th7} this proves $E_\alpha^{0,\rm
sym}(\R^1) \subset F$. It still remains to show that the Gaussian
distribution is in $F$. Since $\mE_\alpha$ is a bijection from
$I(\R^1)$ onto $E_\alpha(\R^1)$ and since $E_\alpha^{0,\rm
sym}(\R^1) \subset F$, it follows that
$$  \{ \mu \in I_{\rm sym} (\R^1): \text{$\mu$ has
Gaussian part $0$}\} = \mE_\alpha^{-1}(E_\alpha^{0,\rm sym} (\R^1))
\subset \mE_\alpha^{-1} (F),$$ where we used Proposition
\ref{prop-1} and Corollary~\ref{Cor-2}. Let $\zeta$ be a standard
normal distribution and consider the compound Poisson distribution
$\mu_n$ with L\'evy density $n\sqrt{\frac{n}{2\pi}} e^{-nx^2/2} \,
dx$. Then $\mu_n$ converges weakly to $\zeta$ as $n\to\infty$ as
shown in the proof of Theorem~8.1(i) of Sato~\cite{S99}, pp.44-45.
Also observe that $\mu_n$ has Gaussian part $0$ and is symmetric, so
that $\mu_n \in \mE_\alpha^{-1} (F)$. Hence  $\mE_\alpha(\mu_n) \in
F$ and $\mE_\alpha(\mu_n)$ converges weakly to $\mE_\alpha(\zeta) =
\law (\sqrt{\Gamma({1+2/\alpha})} \zeta)$, which is in $F$ since $F$
is closed under weak convergence. By the scaling property, also
$\law (\zeta)\in F$, so that $E_\alpha^{\rm sym}(\mathbb R ^1) \subset F$, giving
(iii).

(ii) Let $F$ be a class of infinitely divisible distributions on
$\bR^1$ which is closed under convolution, weak convergence, scaling
and taking of powers, and contains $\law (Z_1^{(\alpha)})$. Then it
follows in complete analogy to the proof of (iii) that
$E_\alpha^{+,0} (\bR^1) \subset F$, now using
Theorem~\ref{theorem:th8}(ii) instead of Theorem~\ref{theorem:th7}.
Since $\law(t^{-1} Z_t^{(\alpha)}) \in F$ for every $t>0$ by the
power and scaling property and since $Z_1^{(\alpha)}$ has finite and
positive expectation, it follows from the strong law of large
numbers that $t^{-1} Z_t^{(\alpha)}$ converges almost surely and
hence in distribution to $E (Z_1^{(\alpha)})>0$ as $t\to \infty$,
implying that $\delta_{E (Z_1^{(\alpha)})} \in F$. Together with the
scaling and convolution property this shows $E_\alpha^+ (\bR^1)
\subset F$, giving (ii).

(i) Let $F$ be a class of infinitely divisible distributions which
is closed under convolution, weak convergence, scaling and taking of
powers and contains $\law (Z_1^{(\alpha)})$ and $\law
(-Z_1^{(\alpha)})$. By (ii) we have $\delta_1 \in F$ and similarly
$\delta_{-1} \in F$, and using Theorem~\ref{theorem:th8}(i) we then
obtain similarly to the proof of (iii) that  $E_\alpha^0 (\bR^1)
\subset F$. Using (iii) we see further see that $F$ must contain the
Gaussian distributions, and it follows that $E_\alpha(\bR^1) \subset
F$, finishing the proof.
\end{proof}

\begin{rem} \label{rem-5.2}
$E_1^+(\bR^1)=B(\bR_+)$ and $E_1(\bR^1)=B(\bR^1)$ are the
Goldie--Steutel--
\linebreak
Bondesson classes on $\bR_+$ and $\bR^1$,
respectively. Both are themselves the smallest classes of infinitely
divisible distributions closed under convolution and weak
convergence and contain the distributions of all (elementary, resp.)
mixtures of exponential variables; see Bondesson~\cite{B82} for
$B(\bR_+)$ and Barndorff-Nielsen et al.~\cite{BNMS06} for
$B(\bR)$. So, Theorems~\ref{theorem:th7} and~\ref{theorem:th8} give
us another interpretation of (subclasses) of $E_\alpha(\bR^1)$ as
(essential limits) of stochastic integrals with respect to some
fixed compound Poisson processes.
Observe that also Theorem~\ref{theorem:th9} gives a new
interpretation of $B(\bR ^1)$ and $B(\bR_+)$, since it is based on
containing the law of some compound Poisson process, which is not an
exponential distribution.
\end{rem}

\begin{rem}
Recently, the authors had a chance to look at the paper by James et
al.~\cite{JRY08}. In their paper, the authors introduced the
Wiener-Gamma integrals, which are stochastic integrals with respect
to the standard gamma process that is a subordinator without drift
with L\'evy measure $x^{-1}e^{-x}{\bf 1}_{(0,\infty)}(x)dx$, and
studied generalized gamma convolutions, which are related to the
Thorin class. Actually, the Thorin class $T(\R ^1)$ in one
dimension, is the smallest class that contains all gamma
distributions and is closed under convolution and weak convergence.
Distributions in the class $T(\R ^1)$ are named generalized gamma
convolutions. As shown in \cite{BNMS06}, $T(\rd)$ is characterized
as $\Psi(I _{\log}(\rd))$, and $B(\rd)$ is characterized as $\Up
(I(\rd))$. In this paper, we characterized $E_1^0(\R ^1)$, as
mentioned in Remark~\ref{rem-5.2}, in a different way by using
compound Poisson processes. From this point of view, the way of
using gamma process in the paper by James et al.~\cite{JRY08} has a
similar fashion.
\end{rem}

\vskip 10mm
\section*{Acknowledgement}
Parts of this paper were written while A.~Lindner was visiting the
Department of Mathematics at Keio University. He thanks for their
hospitality and support.

\end{document}